\documentclass[10pt,twocolumn]{IEEEtran}

\usepackage[T1]{fontenc}
\usepackage{units}
\usepackage{amsthm}
\usepackage{amsmath}
\usepackage{amssymb}
\usepackage{enumerate,color,cite}
\usepackage{graphicx}
\usepackage{caption}
\usepackage[labelformat=simple]{subcaption}
\usepackage{algorithm, algorithmic}
\usepackage[hidelinks]{hyperref}

\graphicspath{{./Figs/}}

\newtheorem{assumption}{Assumption}[section]
\newtheorem{thm}{Theorem}

\newtheorem{cor}{Corollary}
\newtheorem{lem}{Lemma}
\newtheorem{rem}{Remark}

\newcommand{\cS}{\mathcal{S}}
\newcommand{\cN}{\mathcal{N}}
\newcommand{\lpthres}{{t_{\text{thres}}}}
\newcommand{\phyN}[1]{\cN_{#1 \leftarrow}}
\newcommand{\multiN}[1]{\overline{\mathcal{N}}_{#1 \leftarrow}}
\newcommand{\dreachN}[1]{\mathcal{N}_{#1 \rightarrow}}
\newcommand{\reachN}[1]{\overline{\mathcal{N}}_{#1 \rightarrow}}
\newcommand{\commcost}[1]{c_{#1}^{cm}(\multiN{#1}, \{\multiN{l}'\}_{l\in\reachN{#1}})}
\newcommand{\ones}[1]{\mathbf{1}_{#1}}

\begin{document}

\title{Multi-hop Diffusion LMS for Energy-constrained Distributed Estimation}

\author{Wuhua\,Hu,\,\IEEEmembership{Member,~IEEE}, and Wee\,Peng\,Tay,\,\IEEEmembership{Senior Member,~IEEE}%
\thanks{Part of this work was presented at the International Conference on Information Fusion, Jul.\ 2014. This work was supported in part by the Singapore Ministry of Education Academic Research Fund Tier 2 grants MOE2013-T2-2-006 and MOE2014-T2-1-028. W. Hu and W. P. Tay are with the School of Electrical and Electronic Engineering, Nanyang Technological University, Singapore; Emails: \{hwh, wptay\}@ntu.edu.sg}}

\maketitle
\begin{abstract}
We propose a multi-hop diffusion strategy for a sensor network to perform distributed least mean-squares (LMS) estimation under local and network-wide energy constraints. At each iteration of the strategy, each node can combine intermediate parameter estimates from nodes other than its physical neighbors via a multi-hop relay path. We propose a rule to select combination weights for the multi-hop neighbors, which can balance between the transient and the steady-state network mean-square deviations (MSDs). We study two classes of networks: simple networks with a unique transmission path from one node to another, and arbitrary networks utilizing diffusion consultations over at most two hops. We propose a method to optimize each node's information neighborhood subject to local energy budgets and a network-wide energy budget for each diffusion iteration. This optimization requires the network topology, and the noise and data variance profiles of each node, and is performed offline before the diffusion process. In addition, we develop a fully distributed and adaptive algorithm that approximately optimizes the information neighborhood of each node with only local energy budget constraints in the case where diffusion consultations are performed over at most a predefined number of hops. Numerical results suggest that our proposed multi-hop diffusion strategy achieves the same steady-state MSD as the existing one-hop adapt-then-combine diffusion algorithm but with a lower energy budget.
\end{abstract}

\thispagestyle{empty}

\section{Introduction} \label{sec: Introduction}
Distributed estimation arises in a wide range of contexts, including
sensor networks \cite{rabbat2004distributed,kar2009distributed,barbarossa2013distributed}, smart grids \cite{xie2012fully,li2013robust}, machine learning \cite{predd2009collaborative,theodoridis2011adaptive}, and biological networks \cite{passino2002biomimicry,barbarossa2007bio,cattivelli2011modeling}.
Several useful distributed solutions have been developed for this
purpose, such as consensus strategies \cite{tsitsiklis1984convergence,Olfati-Saber2004,kar2009distributed,nedic2010constrained},
incremental strategies \cite{bertsekas1997new,lopes2007incremental},
and diffusion strategies \cite{lopes2008diffusion,cattivelli2010diffusion,sayed2013diffusion,sayed2014diffusion}. The diffusion strategies are particularly attractive because they
are scalable, robust, fully-distributed, and endow networks with real-time
adaptation and learning abilities \cite{sayed2014diffusion}. They
have superior stability ranges and transient performance compared
to the consensus strategies when constant step-sizes are necessary
to enable continuous adaptation under varying network conditions \cite{tu2012diffusion}.
The mean-square stability of diffusion has also been shown to be insensitive to topological
changes caused by asynchronous cooperation among the network nodes \cite{zhao2013asynchronous}.

In each iteration of a diffusion strategy, each node obtains intermediate parameter estimates from its neighbors,\footnote{If a node receives and incorporates an intermediate parameter estimate from another node into its own estimate, we say that the former node consults the latter node.} which are those nodes within communication range of itself. We call these neighboring nodes the \emph{physical neighbors} of the node. The communication cost per iteration of each node in a static network is thus fixed, and the total communication cost can be large if the diffusion algorithm converges slowly. To reduce the number of communication links, \cite{rortveit2010diffusion} and \cite{zhao2012single} limit each node to selecting only one of its neighbors for consultation based on the neighbors' current mean-square deviation (MSD) estimates and a variance-product metric, respectively. Simulations showed that in the steady state these two strategies outperform the probabilistic or gossip alternatives where a single neighbor is randomly selected for consultation at each iteration \cite{boyd2006randomized,lopes2008diffusion_b,takahashi2010link}. The reference \cite{xu2013adaptive} proposed a heuristic algorithm to discount physical neighbors with large excess mean-square errors, while \cite{werner2010energy,arablouei2014distributed} suggested diffusing only a part of the intermediate estimate vector at each iteration so as to reduce the amount of information exchanged, and consequently the communication cost. A game theoretic approach with provable stability was proposed in \cite{namvar2013distributed} for each node to learn in a distributed manner whether to diffuse its estimate based on a utility function that captures the trade-off between its contribution and energy expenditure. A similar idea was also presented in \cite{arroyo2013censoring} with numerical validation.

\begin{figure}
\begin{centering}
\includegraphics[scale=0.65]{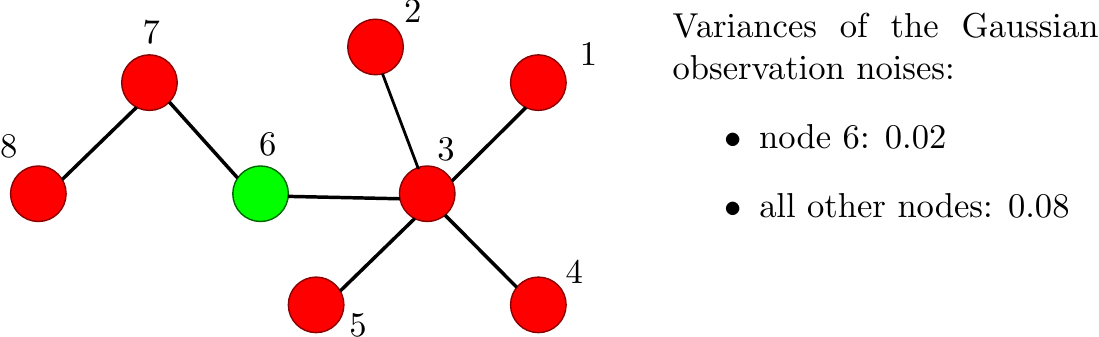}
\par\end{centering}
\caption{A toy example to motivate the use of multi-hop diffusion. Other data used is referred to Section \ref{subsec: numerical-example-tree-network}.}
\label{fig: motivation example}
\end{figure}

{
When designing or upgrading a cooperative sensor network, the strategies in the aforementioned literature are unable to account for predefined node energy budgets even though they are more energy-efficient overall. This prevents efficient energy planning even if we have knowledge about the network working environment (which may be inferred from historical data). Moreover, as these strategies only allow a node to exchange information with its physical neighbors, this limits the estimation performance that a network can achieve.
}
 For example, consider an undirected network, as shown in Fig.~\ref{fig: motivation example}, in which all edges have the same communication length, and the data model, except for the node noise variances, is the same as in the example in Section \ref{subsec: numerical-example-tree-network}. Relying on one-hop communications, the traditional adapt-then-combine (ATC) diffusion strategy \cite{sayed2014diffusion} invokes 8 broadcasts\footnote{
 {
 In this paper, a broadcast means communication of a node with all of its directly reachable neighbors as defined in Section \ref{sec:Problem-formulation}.
 }} per iteration and results in an average steady-state network MSD of -49.0 dB.
In contrast, if we allow two-hop consultations, then the steady-state network MSD can be improved by 2 dB with the same number of broadcasts, or kept the same with only 4 broadcasts by letting nodes 3, 6, and 7 broadcast their intermediate estimates, and node 3 relay the intermediate estimate from node 6 to nodes 1, 2, 4, and 5 at every iteration. In the latter case, the communication cost per iteration is \emph{halved} compared to the ATC diffusion strategy, but the steady-state network MSD remains the same since ``high quality'' information from node 6 is diffused to more nodes than in the ATC strategy.
Although the implementation complexity is somewhat increased (node 3 needs to be programmed to rebroadcast what it receives from node 6), this example shows that to achieve an optimal network MSD-communication cost trade-off requires the use of multi-hop diffusions.

In this paper, we consider diffusion estimation with local and network-wide energy budgets per iteration, and the use of \emph{multi-hop} consultations, i.e., a node that is not a physical neighbor can transmit its intermediate estimate to another node via a relay path in the network within the same iteration step. This is in sharp contrast to all aforementioned literature, which considers only \emph{single-hop} consultations in each iteration. Our main contributions are the following:
\begin{enumerate}[(i)]
	\item We generalize the concept of single-hop diffusion from physical neighbors to multi-hop diffusion from a set of \emph{information neighbors}. In particular, we propose a multi-hop version of the ATC diffusion algorithm, which we call mATC. We formulate and apply mATC to a distributed estimation problem with local and network-wide energy constraints.
	\item For a given set of information neighbors, we provide a rule to select combination weights for mATC that optimizes an approximate trade-off between the convergence rate of the algorithm and the steady-state network MSD.
	\item
	Given the network topology and data and noise variance profiles of each node, we show how to select information neighbors to minimize an upper bound of the steady-state network MSD, subject to local and network-wide energy budgets per iteration, in two network classes: simple networks with a unique transmission path from one node to another, and arbitrary networks utilizing diffusion consultations over at most two hops. We formulate the problem as an offline centralized mixed integer linear program (MILP), and show that the selection is invariant to homogeneous scaling of node observation noise variances.
	\item
	Our MILP requires knowledge of the network topology, and data and noise variance profiles of each node, which may be impractical in some applications. To overcome these requirements, we develop an approximate distributed and adaptive optimization algorithm to select the information neighbors for arbitrary networks utilizing diffusion consultations over at most $h$ hops, in the absence of a network-wide energy budget constraint.
\end{enumerate}

The concept of multi-hop diffusion unifies non-cooperative, distributed diffusion and centralized estimation strategies into a single framework, and allows us to study the trade-offs amongst these strategies easily. Our proposed strategy has the advantage that it achieves good trade-offs between estimation accuracy and \emph{predefined} energy budgets, which the standard diffusion strategy or the approaches in \cite{rortveit2010diffusion,zhao2012single,xu2013adaptive,namvar2013distributed} cannot incorporate. As wireless sensor networks with renewable energy sources become more popular in practical implementations, energy constraints need to be accounted for explicitly in the estimation algorithm \cite{nokleby2013toward}.
{
We also note that multi-hop diffusion is different from geographic gossip \cite{dimakis2008geographic} (or path averaging gossip \cite{benezit2007gossip}), which relies on randomized pair-wise (or relay path-wise) cooperation to achieve average consensus.
}

{
Our proposed multi-hop diffusion strategy requires information relaying and is arguably more complex than the standard diffusion strategy. If each multi-hop relaying is to be completed within each diffusion iteration, as is assumed in our analysis, this may require the nodes to take observations at a slower rate due to a longer communication delay at each iteration. As communication delay is proportional to the number of messages each node needs to relay, which is in turn proportional to the energy budget at each node, our results on the MSD performance-energy trade-off in Section \ref{sec:Numerical-Example} can be interpreted as a trade-off between the MSD performance and delay. In addition, one can perform the relaying over multiple diffusion iterations, so that intermediate estimates of information neighbors more than one hop away are combined only in a later diffusion iteration (we call this \emph{asynchronous} mATC). Our simulation results in Section \ref{sec:Numerical-Example} demonstrate that asynchronous mATC has similar MSD performance as mATC.
}

The rest of this paper is organized as follows. In Section \ref{sec:Problem-formulation},
we introduce our data model and notations, and formulate the energy-constrained distributed optimization problem. In Section \ref{sec: generalized-diffusion-adaptation}, we introduce the concept of multi-hop diffusion adaptation. In Section \ref{sec: selecting the combination weights}, we propose a combination weight to optimize an approximate trade-off between the convergence rate and the steady-state network MSD, and in Section \ref{sec: selecting the info neighbors}, we show how to choose an approximately optimal set of information neighbors for every node in two classes of networks using an offline optimization, and also general networks with only local energy budget constraints in an adaptive procedure. Numerical results and conclusions follow in Sections \ref{sec:Numerical-Example} and \ref{sec: Conclusion}, respectively.

\emph{Notations.} The notation $\mathbb{R}_{\ge0}$ denotes the space of non-negative real numbers, $|\mathcal{N}|$ denotes the cardinality of a discrete set $\mathcal{N}$, $\ones{N}$ represents a vector of size $N$ with all entries equal to one, $I_N$ is an $N\times N$ identity matrix, $A^T$ is the transpose of the matrix $A$, and $\lambda_{max}(A)$ and $\rho(A)$ are the largest eigenvalue and the largest absolute eigenvalue of the matrix $A$, respectively. The operation $A\otimes B$ denotes the Kronecker product of the two matrices $A$ and $B$. The relation $A\succeq (\text{or} \preceq) B$ means that the matrix $A-B$ is positive (or negative) semi-definite, and similarly $A\succ (\text{or} \prec) B$ means the matrix $A-B$ is positive (or negative) definite. The notation $\text{col}\{\cdot\}$ denotes a column vector in which its arguments are stacked on top of each other, $\text{diag}\{\cdot\}$ denotes a diagonal matrix constructed from its arguments. We use boldface letters to denote random quantities (e.g., $\boldsymbol{x}$) and normal letters to denote their realizations or deterministic quantities (e.g., $x$). The symbol $\mathbb{E}\boldsymbol{x}$ denotes the expectation of the random variable $\boldsymbol{x}$, and ``s.t.'' is abbreviation for ``subject to''.

\section{Problem Formulation\label{sec:Problem-formulation}}

We adopt the same notations as in \cite{tu2012diffusion,sayed2014diffusion} for our problem formulation. Consider a network represented by a directed graph $\mathcal{G} = (\mathcal{N},\mathcal{E})$, where $\mathcal{N}=\{1,\,2,\,...,\, N\}$ is the set of nodes, and $\mathcal{E}$ is the set of communication links between nodes.\footnote{An undirected graph is treated as a directed graph by replacing each undirected edge with two edges of opposite directions.} Node $l$ is said to be a \textit{physical neighbor} of node $k$ if $l=k$ or $(l,k)\in\mathcal{E}$, and is said to be within the \emph{multi-hop neighborhood} of node $k$ if there is a path in $\mathcal{G}$ from node $l$ to node $k$. Let the physical neighborhood of node $k$ be $\phyN{k}$, and its multi-hop neighborhood be $\multiN{k}$. We have $\phyN{k}\subseteq\multiN{k}$.

On the other hand, we say that node $l$ is within the \emph{reachable neighborhood} of node $k\neq l$ if there is a path in $\mathcal{G}$ from node $k$ to node $l$. We say that node $l$ is directly reachable from node $k$ if $(k,l)\in\mathcal{E}$. We let $\dreachN{k}$ be the set of directly reachable neighbors of node $k$, and $\reachN{k}$ be the reachable neighborhood of node $k$. We have $\dreachN{k}\subseteq\reachN{k}$. The various types of neighbors are illustrated in Fig.~\ref{fig:different_neighbors}.

At every iteration $i$, each node $k$ is able to observe realizations $\{d_{k}(i),\, u_{k,i}\}$ of a scalar random process $\boldsymbol{d}_{k}(i)$
and a $1\times M$ vector random process $\boldsymbol{u}_{k,i}$ with
a positive definite covariance matrix, $R_{u,k}=\mathbb{E}\boldsymbol{u}_{k,i}^{*}\boldsymbol{u}_{k,i}\succ0$.
The random processes $\{d_{k}(i),\, u_{k,i}\}$ are related via the linear regression model \cite{sayed2014diffusion}:
\[
\boldsymbol{d}_{k}(i)=\boldsymbol{u}_{k,i}\omega^{o}+\boldsymbol{v}_{k}(i),
\]
where $\omega^o$ is an $M\times 1$ parameter to be estimated, and $\boldsymbol{v}_{k}(i)$ is measurement noise with variance $\sigma_{v,k}^{2}$, and assumed to be temporally white and spatially
independent, i.e.,
\[
\mathbb{E}\boldsymbol{v}_{k}^{*}(i)\boldsymbol{v}_{l}(j)=\sigma_{v,k}^{2}\delta_{kl}\delta_{ij},
\]
where $\delta_{kl}$ is the Kronecker delta function. The regression data $\boldsymbol{u}_{k,i}$ are likewise assumed to be temporally white and spatially independent.
The noise $\boldsymbol{v}_{k}(i)$ and the regressors $\boldsymbol{u}_{l,j}$
are assumed to be independent of each other for all $\{k,\, l,\, i,\, j\}$.
All random processes are assumed to be zero mean. The above data model has been frequently used in the parameter estimation literature \cite{sayed2014diffusion}, and are useful in studies of various adaptive filters \cite{sayed2008adaptive}.

The objective of the network is to estimate $\omega^{o}$ in a distributed and iterative way subject
to certain energy constraints. During the iterative estimation process, the energy cost of node $k$
per iteration consists of sensing cost, computing cost and communication
cost (incurred to disseminate or
relay intermediate estimates to the physical neighbors of a
node). While the sensing and computing costs are almost the
same for all nodes, the communication cost depends on the information that is disseminated
or relayed by a node in every iteration and forms the major cost incurred in the estimation process. For simplicity, we ignore the sensing and computing costs and use the terms ``energy cost'' and ``communication cost'' interchangeably throughout the paper. Denote the communication cost per iteration of a node $k$ as $c_{k}^{cm}$. The nodes estimate $\omega^{o}$ by solving a constrained least mean-squares (LMS) problem:
\begin{equation}
\begin{aligned}(\text{P0})\quad &\min_{\omega}\sum_{k\in\mathcal{N}}\mathbb{E}|\boldsymbol{d}_{k}(i)-\boldsymbol{u}_{k,i}\omega|^{2}\\
\text{s.t.,}\,\, &c_{k}^{cm}\le c_{k},\,\,\forall k\in\mathcal{N},\\
&\sum_{k\in\mathcal{N}} c_{k}^{cm} \le c,
\end{aligned}
\label{eq: energy constraints}
\end{equation}
where $c_{k}$ and $c$ are the node- and network-wide energy budgets imposed in each iteration, respectively.

The ATC diffusion strategy solves (P0) without the energy constraints by using the following update equations \cite{sayed2013diffusion,sayed2014diffusion}:
\begin{equation}
\begin{aligned}\boldsymbol{\psi}_{k,i} & =\boldsymbol{\omega}_{k,i-1}+\mu_{k}\boldsymbol{u}_{k,i}^{*}[\boldsymbol{d}_{k}(i)-\boldsymbol{u}_{k,i}\boldsymbol{\omega}_{k,i-1}],\\
\boldsymbol{\omega}_{k,i} & =\sum_{l\in\phyN{k}}a_{lk}\boldsymbol{\psi}_{l,i},
\end{aligned}
\label{eq: original ATC}
\end{equation}
where $\mu_{k}$ is a positive step-size parameter, and $a_{lk}$ are combination weights satisfying
\[
a_{lk}\ge0,\,\,A^{T}\ones{N}=\ones{N},\,\,\text{and }a_{lk}=0\,\,\text{if }l\notin\phyN{k}.
\]
Here, $A=(a_{lk})_{N\times N}$ is the combination weight matrix. The strategy consists of two steps, the \textit{adaptation} step and the \textit{consultation} (also known as the combination or diffusion) step. In the
adaptation step, each node adapts its local estimate to an intermediate estimate $\boldsymbol{\psi}_{k,i}$ by using
the new data available, and the consultation step combines the intermediate
estimates from the physical neighborhood of a node through a weighted
sum to obtain a local estimate $\boldsymbol{\omega}_{k,i}$ for the current iteration. In this paper, we consider only the ATC form of diffusion since it outperforms other alternative diffusion strategies under mild technical conditions  \cite{tu2012diffusion}.

The ATC strategy however does not have the flexibility to take into account the energy constraints in (\ref{eq: energy constraints}). This is because for a given network, the ATC
strategy invokes a fixed communication cost at every node in each
iteration. This motivates us to consider a flexible diffusion strategy, which allows multi-hop consultations under predefined energy budgets.

\section{Multi-hop Diffusion Adaptation \label{sec: generalized-diffusion-adaptation}}

\begin{figure}
\begin{centering}
\includegraphics[scale=0.6]{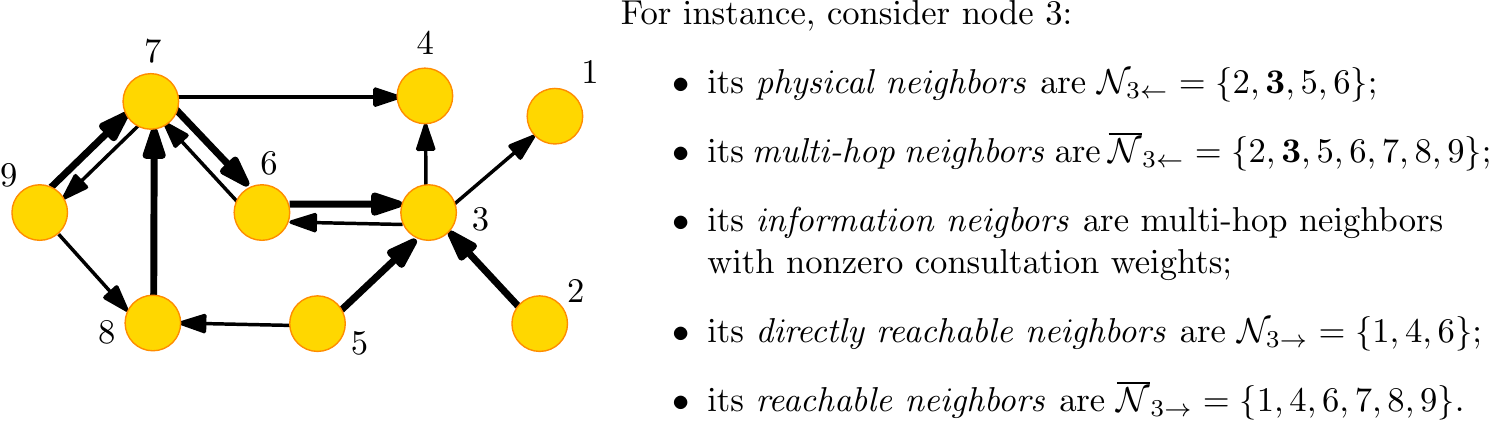}
\par\end{centering}
\caption{The different types of neighbors of a node. \label{fig:different_neighbors}}
\end{figure}

In this section, we extend the ATC strategy by allowing a node to consult any node in its multi-hop neighborhood. The resulting mATC strategy uses the following update equations,
\begin{equation}
\begin{aligned}
\boldsymbol{\psi}_{k,i} & =\mathbf{\omega}_{k,i-1}+\mu_{k}\boldsymbol{u}_{k,i}^{*}[\boldsymbol{d}_{k}(i)-\boldsymbol{u}_{k,i}\boldsymbol{\omega}_{k,i-1}],\\
\boldsymbol{\omega}_{k,i} & =\sum_{l\in\multiN{k}}a_{lk}\boldsymbol{\psi}_{l,i},
\end{aligned}
\label{eq: generalized ATC}
\end{equation}
where the combination weights satisfy
\begin{equation}
a_{lk}\ge0,\,\, A^{T}\ones{N}=\ones{N},\,\,\text{and }a_{lk}=0\,\,\text{if }l\notin\multiN{k}.\label{eq: constraints on combination matrix}
\end{equation}
The only difference between mATC and ATC is in the combination step: the node $k$ consults its multi-hop neighbors $\multiN{k}$, which include the physical neighbors $\phyN{k}$ as a subset. If $a_{lk} > 0$, we say that node $l$ is an \emph{information neighbor} of node $k$ (cf.\ Fig.~\ref{fig:different_neighbors} for an illustration).

This simple modification to the ATC strategy generalizes the diffusion concept to cover centralized estimation at one extreme, and non-cooperative estimation at the other extreme. This unifies the centralized, non-cooperative, and distributed strategies into a single framework, which allows us to study the trade-offs amongst them easily.

The mATC strategy inherits all stability and performance results of
the ATC strategy because the generalization introduced in the combination
matrix $A$ does not affect the analysis. Specifically, the network estimation
is mean stable for any choice of $A$
\textit{if and only if} $\mu_{k}<\frac{2}{\lambda_{\max}(R_{u,k})}$. The same condition holds for mean-square stability if the step sizes $\{\mu_{k}\}_{k\in\mathcal{N}}$ are sufficiently small. Interested readers are referred
to \cite{tu2012diffusion,sayed2013diffusion,sayed2014diffusion}
for proofs of these stability results. Here we only summarize the mean-square
performance results of the mATC strategy as these will be used in the sequel.

Denote the estimation error vector of an arbitrary node $k$ at iteration
$i$ as $\tilde{\boldsymbol{\omega}}_{k,i}\triangleq\omega^{o}-\boldsymbol{\omega}_{k,i}$.
Collect all error vectors and step-sizes across the network into a
block vector and block matrix in
\begin{eqnarray*}
\tilde{\boldsymbol{\omega}}_{i} & \triangleq & \text{col}\{\tilde{\boldsymbol{\omega}}_{1,i},\,\tilde{\boldsymbol{\omega}}_{2,i},\,...,\,\tilde{\boldsymbol{\omega}}_{N,i}\},\\
\mathcal{M} & \triangleq & \text{diag}\{\mu_{1}I_{M},\,\mu_{2}I_{M},\,...,\,\mu_{N}I_{M}\},
\end{eqnarray*}
and let $\mathcal{A}\triangleq A\otimes I_{M}$.
We further define
the block diagonal matrix $\mathcal{R}$ and the $N\times N$ block
matrix $\mathcal{B}$ with blocks of size $M\times M$ each, as follows:
\begin{align*}
\mathcal{R} & \triangleq\mathbb{E}\text{diag}\{\boldsymbol{u}_{1,i}^{*}\boldsymbol{u}_{1,i},\,\boldsymbol{u}_{2,i}^{*}\boldsymbol{u}_{2,i},\,...,\,\boldsymbol{u}_{N,i}^{*}\boldsymbol{u}_{N,i}\},\\
\mathcal{B} & \triangleq\mathcal{A}^{T}(I_{NM}-\mathcal{M}\mathcal{R}).
\end{align*}
Then, the mean network error evolves as $\mathbb{E}\tilde{\boldsymbol{\omega}}_{i}=\mathcal{B}\mathbb{E}\tilde{\boldsymbol{\omega}}_{i-1}$.
For any Hermitian nonnegative-definite weighting matrix $\Sigma$, we have the following approximation up to first order in $\mu_{k}$:
\begin{align}
\mathbb{E}\left\Vert \tilde{\boldsymbol{\omega}}_{i}\right\Vert _{\Sigma}^{2} & \thickapprox\text{Tr}\left(\mathcal{B}^{*}{}^{i+1}\Sigma\mathcal{B}^{i+1}\Omega_{-1}\right)+\sum_{j=0}^{i}\text{Tr}\left(\mathcal{B}^{*j}\Sigma\mathcal{B}^{j}\mathcal{Y}\right)\nonumber \\
 & =\mathbb{E}\left\Vert \tilde{\boldsymbol{\omega}}_{i-1}\right\Vert _{\Sigma}^{2}+\text{Tr}\left(\mathcal{B}^{*i}\Sigma\mathcal{B}^{i}\mathcal{Y}\right)\nonumber \\
 & \,\,\,\,\,-\text{Tr}\left(\left(\mathcal{B}^{*}{}^{i}\Sigma\mathcal{B}^{i}-\mathcal{B}^{*}{}^{i+1}\Sigma\mathcal{B}^{i+1}\right)\Omega_{-1}\right),\label{eq: mean-square error evolution-2}
\end{align}
where $\left\Vert \tilde{\boldsymbol{\omega}}_{i}\right\Vert _{\Sigma}^{2}\triangleq\tilde{\boldsymbol{\omega}}_{i}^{*}\Sigma\tilde{\boldsymbol{\omega}}_{i}$, $\Omega_{-1}\triangleq\mathbb{E}\tilde{\boldsymbol{\omega}}_{-1}\tilde{\boldsymbol{\omega}}_{-1}^{*}$ with $\tilde{\boldsymbol{\omega}}_{-1}$ being the initial estimation error, and
\begin{align*}
\mathcal{Y}&\triangleq\mathcal{A}^{T}\mathcal{MSMA}, \\
\mathcal{S}&\triangleq\text{diag}\{\sigma_{v,1}^{2}R_{u,1},\,\sigma_{v,2}^{2}R_{u,2},\,...,\,\sigma_{v,N}^{2}R_{u,N}\}.
\end{align*}
The recursive relation (\ref{eq: mean-square error evolution-2})
can be used to compute the theoretical transient and steady-state network MSDs.

By specifying $\Sigma$ as $\frac{1}{N}I_{NM}$, the above variance
$\mathbb{E}\left\Vert \tilde{\boldsymbol{\omega}}_{i}\right\Vert _{\Sigma}^{2}$
gives the MSD of the network estimate $\boldsymbol{\omega}_{i}$,
which is an average MSD across the network at the
$i$th iteration, i.e., $\text{MSD}_{i}\triangleq\frac{1}{N}\sum_{k\in\mathcal{N}} \mathbb{E}\left\Vert \tilde{\boldsymbol{\omega}}_{k,i}\right\Vert ^{2}$. In particular, as $i\rightarrow\infty$ the steady-state network MSD
is obtained from (\ref{eq: mean-square error evolution-2}) as
\begin{equation}
\text{MSD}_{\infty}\thickapprox\dfrac{1}{N}\sum_{j=0}^{\infty}\text{Tr}\left(\mathcal{B}^{j}\mathcal{Y}\mathcal{B}^{*j}\right).\label{eq: steady-state network MSD}
\end{equation}
The node and network MSDs are controlled by the quantities $\mathcal{B}$
and $\mathcal{Y}$, both of which are dependent on the combination
matrix $A$. Selection of the combinations weights $a_{l,k}$ can be done in two steps:

\textit{Step 1}: Given an arbitrary set of information neighbors
for each node to consult, we derive analytical forms of the combination
weights that optimize the network performance.

\textit{Step 2}: Given the analytical combination weights derived
in Step 1, we optimize the information neighbor set to be consulted
by each node such that the network MSD is minimized subject to predefined energy budget constraints.

The two steps together determine which nodes are consulted by each node and to what extent it is weighted if a consultation happens. The next section presents a combination rule for determining the weights in Step 1, while we discuss Step 2 in Section \ref{sec: selecting the info neighbors}.

\section{Selecting the Combination Weights\label{sec: selecting the combination weights}}

In this section, we aim to select the combination weight matrix $A$ to optimize the steady-state network MSD in (\ref{eq: steady-state network MSD}), given arbitrary information neighbors of each node in the network. The optimization, however, does not admit an analytical solution and has to be solved numerically in general, which prevents finding an adaptive solution under varying network conditions. To keep the adaptation ability of the network, we make a compromise by seeking for an analytical solution that approximately minimizes an \textit{upper bound} of the steady-state network MSD, given by
\begin{align*}
&\text{MSD}_{\infty}
\thickapprox\dfrac{1}{N}\sum_{j=0}^{\infty}\text{Tr}\left(\mathcal{B}^{j}\mathcal{Y}\mathcal{B}^{*j}\right)
\le\dfrac{\lambda_{\max}(\mathcal{Y})}{N}\sum_{j=0}^{\infty}\text{Tr}\left(\mathcal{B}^{j}\mathcal{B}^{*j}\right)\\
&\le M\lambda_{\max}(\mathcal{Y})\sum_{j=0}^{\infty}\lambda_{\max}^{j}(\mathcal{B}\mathcal{B}^{*})
=\dfrac{M\lambda_{\max}(\mathcal{Y})}{1-\lambda_{\max}(\mathcal{B}\mathcal{B}^{*})}\triangleq\overline{\text{MSD}}_{\infty}.
\end{align*}
The first inequality above uses the positive semi-definiteness
of the matrix $\mathcal{Y}$ and the last equality is due to the fact that $\lambda_{\max}(\mathcal{B}\mathcal{B}^{*})<1$ is necessary to ensure mean and mean-square stability \cite{tu2012diffusion,sayed2013diffusion,sayed2014diffusion}.

Given arbitrary information neighbors for each node, the upper bound $\overline{\text{MSD}}_{\infty}$ can be minimized by minimizing $\beta \geq 0$ subject to the eigenvalue constraint $\lambda_{\max}(\mathcal{Y})/\beta +\lambda_{\max}(\mathcal{B}\mathcal{B}^{*})\le 1$. This however does not admit a closed-form solution. To obtain an explicit solution for the combination matrix $A$, we approximate and decompose the above problem into two subproblems: Firstly, we solve for an approximate solution of $\beta$ by strengthening the eigenvalue constraint into $\lambda_{\max}\left(\mathcal{Y}/\beta+\mathcal{B}\mathcal{B}^{*}\right)\le 1$. We can solve for an approximate solution of $\beta$ using the following semi-definite program (SDP) (refer to Appendix A for the derivation):
{\small
\begin{equation}
\begin{gathered}\beta^{o}= \arg\min_{\beta}\beta\quad\text{s.t.}\\
\left[\begin{array}{cc}
\beta(\mathcal{MSM})^{-1}+\tilde{Q} & \tilde{Q}(\ones{N}\otimes I_{M})\\
(\ones{N}^{T}\otimes I_{M})\tilde{Q} & (\ones{N}^{T}\otimes I_{M})[\tilde{Q}-I_{NM}](\ones{N}\otimes I_{M})
\end{array}\right]\succeq0,
\end{gathered}
\label{eq: SDP for beta}
\end{equation}
}where $\tilde{Q}\triangleq(I_{NM}-\mathcal{MR})^{-2}\succ0$, which
is independent of the combination weight matrix $A$ to be optimized. The SDP is
convex and hence readily solvable by standard SDP solvers.

Secondly, given the solution of $\beta$, we derive an analytical solution of the weight matrix $A$ by minimizing $\mathrm{Tr}(\mathcal{Y})/\beta+\mathrm{Tr}(\mathcal{B}\mathcal{B}^{*})$, which is an upper bound of the original eigenvalue constraint. This is equivalent to the optimization problem in \eqref{eq: relaxed optimization of the network MSD bound}.

\begin{thm}\label{thm: Optimal-combination-rule}
Suppose that $\beta^o$ is the solution of the SDP (\ref{eq: SDP for beta}), and $\alpha^o\triangleq (\beta^o+1)^{-1}$. The combination weights that solve the following optimization problem
\begin{equation}
\begin{gathered}
\min_{A}\alpha^{o}\mathrm{Tr}(\mathcal{Y})+(1-\alpha^{o})\mathrm{Tr}(\mathcal{B}\mathcal{B}^{*})\\
\text{s.t. }a_{lk}\ge0,\,\, A^{T}\ones{N}=\ones{N},\,\,\text{and }a_{lk}=0\,\,\text{if }l\notin\multiN{k}'.
\end{gathered}
\label{eq: relaxed optimization of the network MSD bound}
\end{equation}
are given as follows:
\begin{equation}
a_{lk}=\begin{cases}
\dfrac{\gamma_{l}^{-2}}{\sum_{j\in\multiN{k}'}\gamma_{j}^{-2}}, & \text{if }l\in\multiN{k}',\\
0, & \text{otherwise},
\end{cases}\label{eq: optimized combination weights}
\end{equation}
where the composite variance $\gamma_{l}^{2}$ is defined by
\begin{equation}
\gamma_{l}^{2}\triangleq\alpha^o\mu_{l}^{2}\cdot\sigma_{v,l}^{2}\cdot\mathrm{Tr}(R_{u,l})+(1-\alpha^o)\mathrm{Tr}\left((I_{M}-\mu_{l}R_{u,l})^{2}\right).\label{eq: gamma}
\end{equation}
\end{thm}
\begin{IEEEproof}
Substituting the expression of $\mathcal{Y}$ into the objective function
of problem (\ref{eq: relaxed optimization of the network MSD bound}), we have
\begin{align*}
\lefteqn{\alpha^{o}\mathrm{Tr}(\mathcal{Y})+(1-\alpha^{o})\mathrm{Tr}(\mathcal{B}\mathcal{B}^{*})}\\
 & =\alpha^{o}\sum_{k\in\mathcal{N}}\sum_{l\in\mathcal{N}}a_{lk}^{2}\mu_{l}^{2}\sigma_{v,l}^{2}\mathrm{Tr}(R_{u,l})\\
 & \,\,\,\,+(1-\alpha^{o})\sum_{k\in\mathcal{N}}\sum_{l\in\mathcal{N}}a_{lk}^{2}\mathrm{Tr}\left((I_{M}-\mu_{l}R_{u,l})^{2}\right).
\end{align*}
Therefore problem (\ref{eq: relaxed optimization of the network MSD bound})
can be decoupled into $N$ separate optimization problems of the form:{\small
\[
\begin{gathered}\min_{A}\sum_{l\in\mathcal{N}}a_{lk}^{2}\left[\alpha^{o}\mu_{l}^{2}\sigma_{v,l}^{2}\mathrm{Tr}(R_{u,l})+(1-\alpha^{o})\mathrm{Tr}\left((I_{M}-\mu_{l}R_{u,l})^{2}\right)\right]\\
\text{s.t. }a_{lk}\ge0,\,\, A^{T}\ones{N}=\ones{N},\,\,\text{and }a_{lk}=0\,\,\text{if }l\notin\multiN{k}^{'},
\end{gathered}
\]
}from which \eqref{eq: optimized combination weights} follows, and the proof is complete.
\end{IEEEproof}

We call the closed-form solution of the combination weights given in \eqref{eq: optimized combination weights} as the ``balancing rule'', since it optimizes a trade-off between the diffusion convergence rate (measured through $\mathrm{Tr}(\mathcal{B}\mathcal{B}^{*})$), and the steady-state network MSD (measured through $\mathrm{Tr}(\mathcal{Y})$). This is further explained as follows.

We observe that the steady-state network MSD can be further upper bounded by{
\begin{equation}
\mathrm{MSD}_{\infty}\le\left\{ \begin{array}{c}
\dfrac{M\lambda_{\max}(\mathcal{Y})}{1-r_{1}^{2}}\le\dfrac{M\mathrm{Tr}(\mathcal{Y})}{1-r_{1}^{2}}\triangleq\overline{\text{MSD}}_{\infty}^{a}\\
\dfrac{Mr_{2}}{1-\lambda_{\max}(\mathcal{B}\mathcal{B}^{*})}\le\dfrac{Mr_{2}}{1-\mathrm{Tr}(\mathcal{B}\mathcal{B}^{*})}\triangleq\overline{\text{MSD}}_{\infty}^{b}
\end{array}\right.\label{eq: looser steady-state MSD bound}
\end{equation}
}where the two positive scalars $r_{1}$ and $r_{2}$ are given by
\[
r_{1}\triangleq\rho(I_{NM}-\mathcal{MR}),\,\, r_{2}\triangleq\lambda_{\max}(\mathcal{MSM}).
\]
These two upper bounds can be shown to be minimized by $\mu_{l}^{2}\sigma_{v,l}^{2}\mathrm{Tr}(R_{u,l})$ and $\mathrm{Tr}\left((I_{M}-\mu_{l}R_{u,l})^{2}\right)$, respectively. Therefore the composite variance $\gamma_{l}^{2}$ is a summation of the minimizers of the two looser upper bounds in \eqref{eq: looser steady-state MSD bound}. Observing that $\mathrm{Tr}(\mathcal{B}\mathcal{B}^{*})$ is an upper bound of $\lambda_{\max}(\mathcal{B}\mathcal{B}^{*})$, we can alternatively
interpret the term $\mathrm{Tr}\left((I_{M}-\mu_{l}R_{u,l})^{2}\right)$
as an approximate minimizer of the transient network MSD. The balancing rule can then be interpreted as balancing between minimizing the steady-state and the transient network MSDs, with the balance tuned by varying the coefficient $\alpha^o \in [0, 1]$.

To apply the balancing rule (\ref{eq: optimized combination weights}),
each node $k$ needs to know the composite variances $\{\gamma_{l}^{2}\}_{l\in\multiN{k}'}$,
which depend on the two components $\{\mu_{l}^{2}\sigma_{v,l}^{2}\mathrm{Tr}(R_{u,l}),\,\mathrm{Tr}\left((I_{M}-\mu_{l}R_{u,l})^{2}\right)\}$
of each selected information neighbors. Without knowing them a priori,
each node needs to gather their estimates from its information neighbors and use them to update
the combination rule adaptively. The desired estimates can be obtained as moving averages of their realizations using real-time data (a similar method was used in \cite{zhao2012diffusion} to obtain the adaptive relative-variance rule):
\begin{equation}\label{eq: gamma estimate}
\begin{aligned}
\boldsymbol{\hat{\gamma}}_{l,1}^{2}(i) &= (1-\nu_{l}) \boldsymbol{\hat{\gamma}}_{l,1}^{2}(i-1) + \nu_{l} \left\Vert \boldsymbol{\psi}_{l,i}-\boldsymbol{\omega}_{l,i-1}\right\Vert^{2},\\
\boldsymbol{\hat{R}}_{u,l}(i) &= (1-\nu_{l}) \boldsymbol{\hat{R}}_{u,l}(i-1) + \nu_{l} \boldsymbol{u}_{l,i}^*\boldsymbol{u}_{l,i},\\
\boldsymbol{\hat{\gamma}}_{l,2}^{2}(i) &= \mathrm{Tr}\left((I_{M}-\mu_{l}\boldsymbol{\hat{R}}_{u,l}(i))^{2}\right),\\
\boldsymbol{\hat{\gamma}}_{l}^{2}(i) &= \hat{\alpha}^o \boldsymbol{\hat{\gamma}}_{l,1}^{2}(i) + (1-\hat{\alpha}^o) \boldsymbol{\hat{\gamma}}_{l,2}^{2}(i),
\end{aligned}
\end{equation}
where the symbol $\hat{}$ indicates an adaptive estimate. The quantity $\nu_l\in(0,1)$ is a chosen discount factor, and $\hat{\alpha}^o\in[0,1]$ is a balancing coefficient usually chosen close to one. The balancing coefficient can also be obtained through the SDP in Theorem \ref{thm: Optimal-combination-rule} by replacing the noise variances with empirical estimates. This however requires sending the empirical estimates to a central processor, which can be costly for the network, and is therefore done only infrequently. By the adaptation equation in \eqref{eq: generalized ATC}, we can verify that $\mathbb{E}\left\Vert \boldsymbol{\psi}_{l,i}-\boldsymbol{\omega}_{l,i-1}\right\Vert^{2}=\mu_{l}^{2}\sigma_{v,l}^{2}\mathrm{Tr}(R_{u,l})$, and we arrive at an adaptive implementation of the balancing rule:
\begin{equation} \label{eq: adaptive balancing rule}
\boldsymbol{a}_{lk}(i) = \dfrac{\boldsymbol{\hat{\gamma}}_{l}^{-2}(i)}{\sum_{j\in\multiN{k}'}\boldsymbol{\hat{\gamma}}_{j}^{-2}(i)}.
\end{equation}
Note that the estimates $(\boldsymbol{\psi}_{l}(i),\boldsymbol{\hat{\gamma}}_{l}^{2}(i))$ of each information neighbor are transmitted to node $k$ before it applies the rule.

\begin{rem}
The balancing rule reduces to the relative-variance rule of \cite{zhao2012diffusion} if the coefficient $\alpha^o$ is set to 1. This implies that the relative-variance rule minimizes the bound $\overline{\text{MSD}}_{\infty}^{a}$. The balancing rule
is also related to the two-phase rules proposed in \cite{yu2013a},
where separate combination rules adopted for the transient and
the steady-state phases approximately minimize $\lambda_{\max}(\mathcal{B}\mathcal{B}^{*})$
and $\mathrm{Tr}(\mathcal{Y})$, respectively. In that case, a switching
point between the two-phase rules needs to be estimated online by, e.g., using the technique developed in \cite{Jesus2014adjustment}.
\end{rem}

\section{Selecting the Information Neighbors \label{sec: selecting the info neighbors}}

{
Given the closed-form combination weights derived in the last section,
we now proceed to optimize the information neighbor set of each node
so that the upper bound of the steady-state network MSD, given as $\alpha^o\overline{\text{MSD}}_{\infty}^{a} + (1-\alpha^o)\overline{\text{MSD}}_{\infty}^{b}$, is minimized
under predefined energy constraints.
With the combination weights given in (\ref{eq: optimized combination weights}), the cost function in \eqref{eq: relaxed optimization of the network MSD bound} becomes an explicit function of the composite variances $\{\gamma_{k}\}_{k\in \mathcal{N}}$, which can be further optimized by selecting appropriate information neighbors for each node under the energy budget constraints. This yields, after some algebraic manipulations, the following optimization problem:
\[
\begin{gathered}(\text{P1})\,\,\min_{\{\multiN{k}':\,\multiN{k}'\subseteq\multiN{k}\}_{k\in\mathcal{N}}}\sum_{k\in\mathcal{N}}\dfrac{1}{\sum_{l\in\multiN{k}'}\gamma_{l}^{-2}}\\
\text{s.t. }\commcost{k}\le c_{k},\,\,\forall k\in\mathcal{N},\\
\sum_{k\in\mathcal{N}}\commcost{k} \le c,
\end{gathered}
\]
where $\commcost{k}$ indicates that the communication cost of node $k$ in one iteration depends on its multi-hop neighbors and reachable neighbors since it may be required to relay estimates from its multi-hop neighbors to its reachable neighbors.
}

Problem (P1) is intractable in its current form because of the unknown
information neighbor sets $\{\multiN{k}'\}_{k\in\mathcal{N}}$
and the implicit communication costs $\commcost{k}$ for all $k\in\mathcal{N}$. To obtain a tractable form of (P1), we first introduce binary variables to represent the sets $\{\multiN{k}'\}_{k\in\mathcal{N}}$, and then model the in-network communications to get an explicit expression for $\commcost{k}$. This turns out to be very complex if the network has an arbitrary topology where there are multiple paths from one node to another. To reduce the complexity and make (P1) tractable, we consider two special cases separately:
\begin{itemize}
  \item Case 1 (Simple topology): For any pair of nodes, there is at most one directed simple path connecting them.
	\item Case 2 (Two-hop consultations): The network has an arbitrary topology but the information neighbors of every node are restricted to be within two hops away.
\end{itemize}

In the sequel, we derive an explicit form of (P1) for each of the two cases as an MILP, which is solved offline before running mATC on the network. This applies in the case where there is a decision center with knowledge of the entire sensor network topology and has inferred the data and noise profiles of each node using historical data. To overcome these requirements, we also show that the general problem (P1) admits an approximate and distributed solution that can be obtained online if only local energy budget is imposed on each node.

\subsection{Explicit formulation as an MILP}

In this subsection, we introduce binary and auxiliary
continuous decision variables to reformulate problem (P1) into an MILP that is solvable using standard solvers.

Before reformulating (P1), we first derive an explicit expression for the communication cost $\commcost{k}$ by introducing two classes of binary variables. The first
class of binary variables are the selection variables $\delta_{lk}\in\{0,\,1\}$,
for all $l\in\multiN{k}$ and $k\in\mathcal{N}$, where $\delta_{lk}=1$ if and only if node $l$ is selected to be an information neighbor of node $k$. We note that $a_{lk}=0$ if and only if $\delta_{lk}=0$. The second class of binary variables are
the relay variables $\pi_{lk}\in\{0,\,1\}$, for all $l\in\multiN{k}$
and $k\in\mathcal{N}$, where $\pi_{lk}=1$ if and only if
node $k$ relays the information originating from node $l$. We have $\pi_{kk}=1$ if and only if node $k$ broadcasts its own intermediate estimate. We make the following assumptions.

\begin{assumption}\label{assumpt:broadcast}
Every broadcast conveys information from a single node, and incurs a communication
cost (which may be different for different nodes). All nodes having the broadcast node as a physical neighbor receives the information being broadcast.
\end{assumption}

\begin{assumption}\label{assumpt:relay}
At every iteration, each node relays the same piece of information at most once.
\end{assumption}

With the above two assumptions, the communication cost, $\commcost{k}$, of node $k$ in a single iteration is equal to the number of intermediate estimates it needs to relay and diffuse, multiplied by the energy cost incurred in each broadcast, i.e.,
\begin{equation}\label{eqn:commcost}
\commcost{k}=c_{k}^{cm,0}\sum_{l\in\multiN{k}}\pi_{lk},\forall k\in\mathcal{N},
\end{equation}
where $c_{k}^{cm,0}$ is the constant cost per broadcast by node $k$. Note that $\pi_{lk}$ is a function of $\{\multiN{j}': j\in\reachN{k}\cap\reachN{l}\}$.

We now investigate the relationships amongst the relay variables $\pi_{lk}$
and the selection variables $\delta_{lj}$ and then reformulate (P1) into an MILP, under the two special cases alluded to above.

\subsubsection{Case 1 (Simple topology)}

In this case, the unique directed path from a node $l$ to a reachable neighbor $j$ is characterized by a set of binary constants, $\{\eta_{lj,k}\}_{k\in\multiN{j}\backslash\{j\}}$, where $\eta_{lj,k}=1$ if and only if node $k$ is on the path from node $l$ to node $j$. Then, the selection variables $\delta_{lj}$ and the relay variables $\pi_{lk}$ are related to each other as follows:
\begin{equation}
\pi_{lk}=\min\left\{ 1,\,\,\sum_{j\in\reachN{l}\backslash \{k\}}\eta_{lj,k}\delta_{lj}\right\},\forall l\in\multiN{k},k\in\mathcal{N},\label{eq: simple-topology-pi-delta-a}
\end{equation}
which implies that node $k$ relays node $l$'s information if and only if node $l$ is consulted by some node $j$, and node $k$ ($\neq j$) is on the directed path from node $l$ to node $j$. The bound by 1 in the relation \eqref{eq: simple-topology-pi-delta-a} is due to Assumption \ref{assumpt:relay}, in which we assume that node $l$'s estimate is relayed at most once by node $k$ in each iteration. The above
relation \eqref{eq: simple-topology-pi-delta-a} can be further written in the linear form (\ref{eq: cons-simple-topology-pi-delta}) ahead.

The objective function of (P1) can be transformed into a linear function by introducing a couple of linear constraints. Using the selection variables $\delta_{lk}$, we express
the objective function equivalently as $\sum_{k\in\mathcal{N}}\left(\sum_{l\in\multiN{k}}\delta_{lk}\gamma_{l}^{-2}\right)^{-1}$, where the candidate information neighbor set $\multiN{k}'$
is replaced by the multi-hop neighbor set $\multiN{k}$ with the help of selection variables. We then introduce auxiliary optimization variables
\begin{align}\label{aux_z_p}
&z_k = \left(\sum_{l\in\multiN{k}}\delta_{lk}\gamma_{l}^{-2}\right)^{-1}, \textrm{ and }
p_{lk} = z_k \delta_{lk},
\end{align}
which allows us to rewrite the objective as $\sum_{k\in\mathcal{N}}z_{k}$, and to perform McCormick linearization \cite{mccormick1976computability} on the bilinear constraints $p_{lk} = z_k \delta_{lk}$ without loss of optimality.

Consequently, (P1) is equivalent to an MILP problem defined in
(\ref{eq: P2-start})-(\ref{eq: P2-end}), which is called (P2) hereafter.
The data and variables of problem (P2) are referred to Table \ref{tb: symbols for (P2) and (P3)}, which also contains those used in the later problem (P3) for Case 2. Constraints (\ref{eq: cons-obj-a})-(\ref{eq: cons-obj-e})
arise from the linearization of the nonlinear objective of (P1). Constraints (\ref{eq: cons-simple-topology-pi-delta})
describe the relation between the relay and the selection
variables. Constraints (\ref{eq: cons-node-wide-energy}) and (\ref{eq: cons-network-wide-energy})
characterize the energy costs and their budgets for the distributed
estimation in a single iteration.

\begin{table}[!t]
\caption{Symbols used in formulating (P2) and (P3)}
\centering{}%
\label{tb: symbols for (P2) and (P3)}
\begin{tabular}{ll}
\hline
\multicolumn{2}{l}{\textbf{Problem data}}  \tabularnewline
$\mathcal{N}$ & full node set\tabularnewline
$\phyN{k}$ & the set of all physical (i.e., one-hop) neighbors of node $k$\tabularnewline
$\multiN{k}$ & the set of all multi-hop neighbors of node $k$\tabularnewline
$\multiN{k}^2$ & the set of all two-hop neighbors of node $k$\tabularnewline
$\dreachN{k}$ & the set of all directly reachable neighbors of node $k$\tabularnewline
$\reachN{k}$ & the set of all reachable neighbors of node $k$\tabularnewline
$\mathcal{P}_k$ & the set of relay customers of node $k$\tabularnewline
$\mathcal{P}^k$ & the set of relay servers of node $k$\tabularnewline
$\sigma_{v,k}$ & the variance of measurement noise, equal to
 $\mathbb{E}\boldsymbol{v}^*_{k}(i)\boldsymbol{v}_{k}(i)$\tabularnewline
$R_{u,k}$ & the input covariance matrix, equal to
 $\mathbb{E}\boldsymbol{u}^*_{k,i}\boldsymbol{u}_{k,i}$\tabularnewline
 $\mu_k$ & the step size parameter of node $k$ implementing mATC in \eqref{eq: generalized ATC}\tabularnewline
$\gamma_k$ & the composite variance of node $k$ defined in \eqref{eq: gamma} \tabularnewline
$\eta_{lj,k}$ & binary indicator of whether node $k$ is on the transmission\tabularnewline
& path from node $l$ to node $j$ \tabularnewline
$\underline{z}_{k}$ & lower bound of $z_k$ computed as $(\sum_{l\in\multiN{k}}\gamma_{l}^{-2})^{-1}$\tabularnewline
$\bar{z}_{k}$ & upper bound of $z_k$ computed as $(\min_{l\in\multiN{k}}\gamma_{l}^{-2})^{-1}$\tabularnewline
 $c_{k}^{cm,0}$ & the constant energy cost per broadcast by node $k$\tabularnewline
 $c_{k}$ & the energy budget available per iteration for node $k$\tabularnewline
 $c$ & the energy budget available per iteration for the network\tabularnewline
& \tabularnewline
\multicolumn{2}{l}{\textbf{Problem variables}}  \tabularnewline
$\delta_{lk}$ & binary indicator of whether $l$ is consulted by node $k$\tabularnewline
$\pi_{lk}$ & binary indicator of whether $l$'s info is relayed by node $k$\tabularnewline
$p_{lk}$, $z_{k}$ & real auxiliary variables, see \eqref{aux_z_p}\tabularnewline
\hline
\end{tabular}
\end{table}

\begin{figure}[!ht]
\begin{align}
(\text{P2}) & \,\,\min\sum_{k\in\mathcal{N}}z_{k}\label{eq: P2-start}\\
\text{s.t.} & \sum_{l\in\multiN{k}}\gamma_{l}^{-2}p_{lk}=1,\,\,\forall k\in\mathcal{N}\label{eq: cons-obj-a}\\
 & p_{lk}\le\bar{z}_{k}\delta_{lk},\,\,\forall l\in\multiN{k},\, k\in\mathcal{N}\\
 & p_{lk}\ge\underline{z}_{k}\delta_{lk},\,\,\forall l\in\multiN{k},\, k\in\mathcal{N}\\
 & p_{lk}\ge z_{k}+\bar{z}_{k}(\delta_{lk}-1),\,\,\forall l\in\multiN{k},\, k\in\mathcal{N}\\
 & p_{lk}\le z_{k}+\underline{z}_{k}(\delta_{lk}-1),\,\,\forall l\in\multiN{k},\, k\in\mathcal{N}\label{eq: cons-obj-e}\\
 & \pi_{lk}\le\hspace{-8pt}\sum_{j\in\reachN{l}\backslash \{k\}}\hspace{-8pt}\eta_{lj,k}\delta_{lj}\le|\reachN{l}|\pi_{lk}, \,\,\forall l\in\multiN{k},\, k\in\mathcal{N}\label{eq: cons-simple-topology-pi-delta}\\
 & \sum_{l\in\multiN{k}} c_{k}^{cm,0}\pi_{lk}\le c_{k},\,\,\forall k\in\mathcal{N}\label{eq: cons-node-wide-energy}\\
 & \sum_{k\in\mathcal{N}} \sum_{l\in\multiN{k}} c_{k}^{cm,0}\pi_{lk}\le c\label{eq: cons-network-wide-energy}\\
 & \delta_{lk}\in\{0,\,1\},\,\,\forall l\in\multiN{k},\, k\in\mathcal{N}\label{eq: cons-decision-variable-a}\\
 & \pi_{lk}\in\{0,\,1\},\,\,\forall l\in\multiN{k},\, k\in\mathcal{N}\\
 & p_{lk}\in\mathbb{R}_{\ge0},\,\,\forall l\in\multiN{k},\, k\in\mathcal{N}\\
 & z_{k}\in\mathbb{R}_{\ge0},\,\,\forall k\in\mathcal{N}.\label{eq: P2-end}
\end{align}
\end{figure}

Note that problem (P2) is solved offline before the mATC diffusion procedure. It requires a centralized processor to have prior knowledge of the network topology, and data and noise variance profiles of every node, which restricts the frequency that this optimization can be performed. This restriction is however alleviated to some extent by the following result.

\begin{lem}\label{lm: solution-invariance}
The optimal
solution $\{(\delta_{lk},\,\pi_{lk}),\,\,\forall l\in\multiN{k},\, k\in\mathcal{N}\}$
for (P2) is invariant to a homogeneous scaling of the composite variances $\gamma_{l}$ for all $l\in\mathcal{N}$.
\end{lem}
\begin{IEEEproof}
Let the optimal solution of (P2) be $\{(\delta_{lk}^{\star},\,\pi_{lk}^{\star},\, p_{lk}^{\star},\, z_{k}^{\star}),\,\,\forall l\in\multiN{k},\, k\in\mathcal{N}\}$. If $\gamma_{l}$, for every $l\in\mathcal{N}$, is scaled by a constant $\alpha$ to $\alpha\gamma_{l}$,
then it is easy to verify that the optimal solution becomes $\{(\delta_{lk}^{\star},\,\pi_{lk}^{\star},\,\alpha^{2}p_{lk}^{\star},\,\alpha^{2}z_{k}^{\star}),\,\,\forall l\in\multiN{k},\, k\in\mathcal{N}\}$, so that $\{(\delta_{lk},\,\pi_{lk}),\,\,\forall l\in\multiN{k},\, k\in\mathcal{N}\}$ remains unchanged. The lemma is now proved.
\end{IEEEproof}

Lemma \ref{lm: solution-invariance} shows that if changes in the data and noise variances happen
uniformly over all network nodes, then the optimal information neighbor configuration remains unchanged.
In general, the offline or centralized optimization in (P2) is performed infrequently, and is based on historical estimates of the data and noise variance profiles maintained by every node in the network.

\subsubsection{Case 2 (Two-hop consultations)}

In this case, we do not make any assumptions about the network topology but restrict the selection of information neighbors to amongst the multi-hop neighbors that are at most two hops away from each node. We denote the set of neighbors of node $k$ within two hops away as $\multiN{k}^{2}$.

\begin{figure}
\begin{centering}
\includegraphics[scale=0.7]{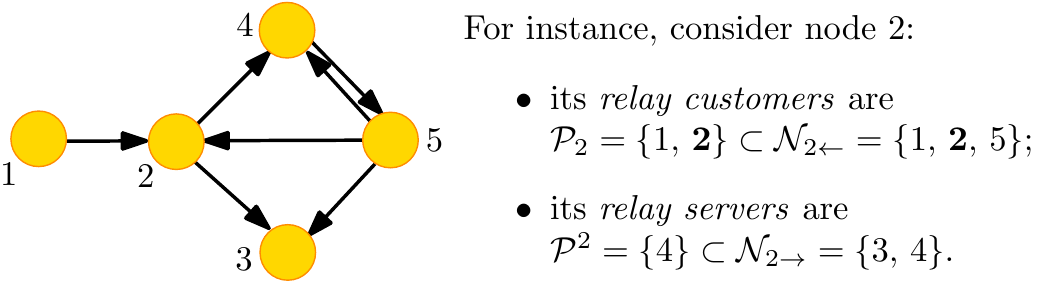}
\par\end{centering}
\caption{The relay customer and server sets of a node. Node 5 is not a relay customer of node 2 because it can reach all its reachable neighbors without node 2's help. Node 3 is not a relay server of node 2 because none of its reachable neighbors needs node 3's help to obtain information from node 2.  \label{fig:relayer_set}}
\end{figure}

With the above restriction, each node is only able to relay information originating from its physical neighbors. We use $\mathcal{P}_{k}$ to denote the set of physical neighbors of node $k$ that have reachable neighbors who may need node $k$ to relay information to.
We call these the \emph{relay customers} of node $k$. On the other hand, suppose that $k \to l \to m$ is a directed path so that node $m$ is reachable but not directly reachable from node $k$. Then, we say that node $l$ is a \emph{relay server} of node $k$.
 Let $\mathcal{P}^{k}$ denote the set of relay servers of node $k$. We have
\begin{align*}
\mathcal{P}_{k} & =\{k\}\cup\{l\in\phyN{k}:\,\dreachN{k}\nsubseteq\dreachN{l}\}\subseteq\phyN{k},\\
\mathcal{P}^{k} & =\{l\in\dreachN{k}:\,\dreachN{l}\nsubseteq\dreachN{k}\}\subseteq\dreachN{k},
\end{align*}
which are illustrated in Fig. \ref{fig:relayer_set}. Then, it is sufficient to define the relay variables $\pi_{lk}\in\{0,\,1\}$ for all $l\in\mathcal{P}_{k},\, k\in\mathcal{N}$, and the selection variables $\delta_{lk}\in\{0,\,1\}$ for all $l\in\multiN{k}^{2},\, k\in\mathcal{N}$.

Next we model the relations between these two sets of variables. We note that, if node
$k$ consults node $l\neq k$, then firstly node $l$ must broadcast its information and secondly,
at least one of its relay servers must relay the information to node $k$ if it
is two hops away from node $k$. These observations are represented mathematically by
\begin{align}
\delta_{lk} & \le \pi_{ll},\,\,\forall l\in\multiN{k}^{2}\backslash\{k\},\, k\in\mathcal{N},\label{eq: arbitrary-topology-delta-pi-a}\\
\delta_{lk} & \le \sum_{j\in\mathcal{P}^{l}\cap\phyN{k}\backslash \{k\}}\pi_{lj},\,\,\forall l\in\multiN{k}^{2}\backslash\phyN{k},\, k\in\mathcal{N}.\label{eq: arbitrary-topology-delta-pi-b}
\end{align}
The first inequality excludes $l=k$ because the inequality
does not hold in that case. The second inequality is defined for nodes $l$ and $k$ that are two hops away from each other, because it is trivially true due to the first inequality if the two nodes are one hop away from each other. These two sets of inequalities give a complete description
of the relations between the relay variables and the selection variables.

Consequently, the information neighbor selection problem (P1) can again be reformulated into an MILP in the form of (P2) with some changes: the constraints (\ref{eq: cons-simple-topology-pi-delta}) are replaced with the new ones in (\ref{eq: arbitrary-topology-delta-pi-a})-(\ref{eq: arbitrary-topology-delta-pi-b}),
the information neighbor set $\multiN{k}$ is replaced with $\multiN{k}^{2}$ everywhere, and furthermore, the relay
variables $\pi_{lk}$ are defined only for all $l\in\mathcal{P}_{k},\, k\in\mathcal{N}$.
We name the new formulation as (P3), to distinguish it from (P2).

Lemma \ref{lm: solution-invariance} similarly holds for the MILP (P3). We also remark that (P3) is applicable to a network having the simple topology assumed in (P2), in which case the constraints
\eqref{eq: cons-simple-topology-pi-delta} (applied to the two-hop neighborhoods) can be treated as valid
inequalities for (P3) to enhance the solution process.

\subsubsection{Valid inequalities to enhance the MILPs}
As combinatorial problems, problems (P2) and (P3) are NP-hard in general.
A typical way to speed up the solution process is by exploiting valid
inequalities (i.e., constraints that maintain the optimal solution), which reduce the search space and enhance the branch-and-cut algorithm used by MILP solvers \cite{golden2008vehicle,burer2012milp}. The following property of (P2) and (P3) at optimality is useful for deriving such valid inequalities.

\begin{lem}\label{lm: monototicity of (P2)-(P3)}
Given a feasible solution of (P2) or (P3), a better solution can be obtained if more information neighbors can be included without violating the energy budgets \eqref{eq: cons-node-wide-energy}-\eqref{eq: cons-network-wide-energy}.
\end{lem}
\begin{IEEEproof}
The conclusion is clear from the objective function of (P1), which is decreasing
in the size of the information neighbor set for each node. Since the objectives of (P1) and (P2) (or (P3)) are equivalent, the conclusion follows immediately.
\end{IEEEproof}


\begin{cor}\label{Cor:valid inequalities}
The following equalities and inequalities are valid for (P2) or (P3):
\begin{align}
& \text{(P2) and (P3): } \delta_{kk}=1,\,\,\forall k\in\mathcal{N},\label{eq: vc-delta-kk}\\
& \text{(P2): } \pi_{lk}\le\delta_{lj},\,\,\forall l\in\multiN{k}\backslash \{j\},\, k\in\phyN{j},\, j\in\mathcal{N},\label{eq: vc-delta-implied-by-pi-P2}\\
& \text{(P3): } \pi_{lk}\le\delta_{lj},\,\,\forall l\in\mathcal{P}_{k}\backslash\{j\},\, k\in\phyN{j},\, j\in\mathcal{N}.\label{eq: vc-delta-implied-by-pi-P3}
\end{align}
\end{cor}
\begin{IEEEproof}
The equalities \eqref{eq: vc-delta-kk} are because for optimality, a node always consults itself, which is evident from Lemma \ref{lm: monototicity of (P2)-(P3)}. The inequalities \eqref{eq: vc-delta-implied-by-pi-P2} (or \eqref{eq: vc-delta-implied-by-pi-P3}) mean if a node relays a piece of information, then every node that has the node as a physical neighbor must use the information in order to achieve optimality. This is proved as follows. By Assumption \ref{assumpt:broadcast}, if a node relays a piece of information, then every node that has the node as a physical neighbor
will receive the information; then by Lemma \ref{lm: monototicity of (P2)-(P3)},
every node receiving the information must use it to ensure optimality. This completes the proof.
\end{IEEEproof}

An additional set of valid inequalities holds for (P2), which states that if node $k$ relays information originating from node $l$ ($\neq k$),
then any predecessor of node $k$ on the unique transmission path must relay the same information towards $k$. Mathematically this means the following:
\begin{lem}\label{lm: additional valid ineq}
The following inequalities hold for (P2):
\begin{align}
\pi_{lk}\le\pi_{lj},\,\,&\forall j\in\{i\in\multiN{k}: \eta_{lk,i}=1\},\,l\in\multiN{k}\backslash\{k\},\,k\in\mathcal{N}.\label{eq: vc-pi-pi-P2}
\end{align}
\end{lem}
\begin{IEEEproof}
We prove the lemma by contradiction. Let there be a pair of nodes $(k,j)$ in the defined set with $1=\pi_{lk}>\pi_{lj}=0$. This implies that there is a simple path that transmits the information from node $l$ to node $k$ while not traversing node $j$. Since by the definition, there is another path that transmits the same information to node $k$ while traversing node $j$, this contradicts with the topological condition in Case 1 that there is a unique transmission path from node $l$ to node $k$. This completes the proof.
\end{IEEEproof}

Within the two-hop neighborhoods, similar inequalities $\pi_{lk}\le \pi_{ll},\,\forall l\in\mathcal{P}_{k}\backslash\{k\},\, k\in\mathcal{N}$, are valid for (P3). They are, however, implied by the valid inequalities (\ref{eq: vc-delta-implied-by-pi-P3}) and the model inequalities (\ref{eq: arbitrary-topology-delta-pi-a})-(\ref{eq: arbitrary-topology-delta-pi-b}), and hence not treated as independent valid inequalities for (P3).

The valid inequalities (\ref{eq: vc-delta-kk})-(\ref{eq: vc-pi-pi-P2})
are added to (P2) or (P3) to reduce the search space and hence speed up the solution
process. However, the computational time may still be prohibitively long if the network has a large size and each node has many reachable neighbors. In that case, we provide in the following subsection, approximate solutions for (P2) and (P3) that can be obtained efficiently by solving linear programs.

\subsection{Approximate solutions via relaxations\label{sub: approximate solution}}

In this subsection, we present a linear programming (LP) approach to solving the MILPs (P2) and (P3) approximately. A typical way of finding an approximate solution for (P2) or (P3) (including all valid inequalities in Corollary \ref{Cor:valid inequalities} and Lemma \ref{lm: additional valid ineq})\footnote{It is useful to include all valid inequalities into the relaxed problem, because they maintain the structural information of the original problem which may otherwise be lost after relaxation.} is by relaxing the binary variables
as continuous variables taking values in $[0,1]$, and then solving the resulting LP. The continuous solution is then translated back into binaries via thresholding. The choice of the threshold needs to ensure that the binary solution obtained satisfies the \emph{hard} energy budget constraints. In the following, we propose a procedure to determine the threshold iteratively for both (P2) and (P3).

\begin{algorithm}[!bt]
\caption{An approximate LP solution for (P2) or (P3)}\label{alg: appro_solution}
\begin{algorithmic}[1]
\STATE{For each $k\in\mathcal{N}$, let $\cS_k = \multiN{k}$ and $\cS_k = \multiN{k}^2$ for (P2) and (P3), respectively.}
\STATE{Relax all binary variables in (P2) or (P3) to continuous variables in $[0,1]$. Solve the resulting LP to obtain $\{(\tilde{\pi}_{lk},\tilde{\delta}_{lk}) : l\in\cS_{k}, k\in\cN\}$ as the optimal relay and selection variables for the relaxed LP.}
\STATE{Initialize the threshold $\lpthres = 1$, the threshold decrement size $\epsilon = \min \{|\tilde{\pi}_{ij} - \tilde{\pi}_{i'j'}| : i\in\cS_{j}, i'\in\cS_{j'}, j,j'\in\cN\}$, and set $\mathrm{EXIT_{flag}} = 0$.}
\WHILE{$\mathrm{EXIT_{flag}} = 0$}\label{start_reduce}
\STATE{For all $l\in\cS_k$ and $k\in\cN$, let $\pi_{lk} = 1$ if $\tilde{\pi}_{lk} \geq \lpthres$ and $\pi_{lk} = 0$ otherwise.} \label{alg-steps: initial-regularized-solution}
\WHILE{network-wide energy budget \eqref{eq: cons-network-wide-energy} not satisfied} \label{alg-steps: global-budget-satisfication-start}
\STATE{Find $(l^*,k^*) = \arg\min_{l,k} \{\tilde{\pi}_{lk} : \pi_{lk} = 1\}$. Set $\pi_{l^*k^*}=0$ and $\mathrm{EXIT_{flag}} = 1$.}
\ENDWHILE \label{alg-steps: global-budget-satisfication-end}
\STATE{Update the threshold as $\lpthres \leftarrow \lpthres - \epsilon$.}
\ENDWHILE\label{end_reduce}
\FOR{each $k\in\cN$}\label{start_node_energy}
\IF{local energy budget \eqref{eq: cons-node-wide-energy} not satisfied for node $k$}  \label{alg-steps: local-budget-satisfication-start}
\STATE{Find $l^* = \arg\min_{l} \{\tilde{\pi}_{lk} : \pi_{lk} = 1\}$. Set $\pi_{l^*k} = 0$.}
\STATE{For (P2), set $\pi_{l^*j} = 0$, for all nodes $j$ reachable from node $l^*$ via node $k$. For (P3), if $l^*=k$, set $\pi_{l^*j} = 0$, for all nodes $j\in\dreachN{k}$.}
\ENDIF \label{alg-steps: local-budget-satisfication-end}
\ENDFOR \label{end_node_energy}
\STATE\label{alg-steps: solution-of-selection-variables}{Determine values of $\{\delta_{lk}\}_{l\in\cS_{k}, k\in\cN}$ from $\{\pi_{lk}\}_{l\in\cS_{k}, k\in\cN}$.}%
\STATE{Return $\{({\pi}_{lk},{\delta}_{lk}): l\in\cS_{k}, k\in\cN\}$ as the approximately optimal solution of the binary relay and selection variables for the original MILP.}
\end{algorithmic}
\end{algorithm}

A common procedure to translate the solutions of relaxed (P2) and (P3) is summarized in Algorithm \ref{alg: appro_solution}. The procedure decreases the threshold gradually (from the maximal value of 1) in lines \ref{start_reduce} to \ref{end_reduce} until it cannot be smaller without violating the network-wide energy budget. Then, in lines \ref{start_node_energy} to \ref{end_node_energy}, the relay variables are adjusted so that local node energy budgets are satisfied. The validity of Algorithm \ref{alg: appro_solution} is proven in the following theorem.

\begin{thm}\label{thm: validity of Algorithm 1}
Algorithm \ref{alg: appro_solution} returns a feasible solution for the relay and selection variables of (P2) and (P3).
\end{thm}
\begin{IEEEproof}
See Appendix A.
\end{IEEEproof}

\begin{rem}
The optimization formulation (P2) or (P3) provides a flexible platform to include additional constraints that may appear in various applications. For example, in (P2), we can easily restrict the information neighbors of a particular node to be within a given number of hops. For another example, the mATC strategy determined by (P3) can lead to improvements over various modified ATC strategies \cite{rortveit2010diffusion,zhao2012single,sayed2013diffusion} by simply imposing constraints on the hops for gathering consultations or on the number of active links associated with each node.
\end{rem}

\subsection{Heuristic distributed solution satisfying local energy budgets \label{sub: distributed solution}}
In this subsection, we present an online distributed algorithm to heuristically solve problem (P1) without the network-wide energy budget, for an arbitrary network in which consultations are restricted to at most $h$ hops. In this case, by Lemma \ref{lm: monototicity of (P2)-(P3)} every node will try to use up its energy budget in each iteration to diffuse intermediate estimates, such that the global cost $\sum_{k\in\mathcal{N}}(\sum_{l\in\multiN{k}^{h'}}\gamma_{l}^{-2})^{-1}$ is minimized, where $\multiN{k}^{h'}$ is a subset of multi-hop neighbors at most $h$ hops away from node $k$. Since each node $k$ contributes to a cost that is decreasing in the size of its information neighborhood $\multiN{k}^{h'}$ and increasing in the composite variance $\gamma_{l}^2$ of each information neighbor $l$, heuristically, we see that each node should broadcast its current intermediate estimate and empirical estimate of its composite variance at each iteration, as well as information from its $h$-hop neighborhood corresponding to nodes with the smallest empirical composite variances, up to its local energy budget. The distributed and adaptive mATC diffusion algorithm is described formally in Algorithm \ref{alg: distributed_solution}.

\begin{algorithm}[!ht]
\caption{Approximate, distributed and adaptive solution for (P1) with only local energy budget constraints}\label{alg: distributed_solution}
\begin{algorithmic}[1]
\STATE{Each node $k\in\mathcal{N}$ initializes the estimates $\boldsymbol{\hat{\gamma}}_{k,1}(-1)$, $\boldsymbol{w}_{k,-1}$ and $\boldsymbol{\hat{R}}_{u,k}(-1)$ (cf. \eqref{eq: gamma estimate}) as a zero scalar, $M\times 1$ vector and $M\times M$ matrix, respectively, and sets $i=0$.}

\STATE{Each node $k\in\mathcal{N}$ uses its own measurement $(\boldsymbol{u}_{k,i}, \boldsymbol{d}_{k}(i))$ to update $\boldsymbol{\psi}_{k,i}$ and $\boldsymbol{\hat{\gamma}}_{k}^{2}(i)$ by \eqref{eq: generalized ATC} and \eqref{eq: gamma estimate}, respectively.} \label{alg-steps: obtain-intermediate-estimate}

\STATE{If its energy budget permits, each node $k\in\mathcal{N}$ broadcasts its estimates $(\boldsymbol{\psi}_{k,i}, \boldsymbol{\hat{\gamma}}_{k}^{2}(i))$ to its directly reachable neighbors, and meanwhile waits for a predefined period of time proportional to $h$ to receive $(\boldsymbol{\psi}_{l,i}, \boldsymbol{\hat{\gamma}}_{l}^{2}(i))$ from its multi-hop neighbors at most $h$ hops away. As allowed by its energy budget, node $k$ then chooses to rebroadcast $(\boldsymbol{\psi}_{l,i}, \boldsymbol{\hat{\gamma}}_{l}^{2}(i))$ if $\boldsymbol{\hat{\gamma}}_{l}^{2}(i-1)$ is among the smallest empirical composite variances it received in the $(i-1)$-th iteration.}

\STATE{Each node $k\in\mathcal{N}$ uses the empirical composite variances received in the previous step to compute the combination weights using \eqref{eq: adaptive balancing rule}, and then combines the received intermediate estimates by \eqref{eq: generalized ATC} to get its own estimate $\boldsymbol{w}_{k,i}$.}

\STATE{Each node completes the $i$-th iteration, sets $i\leftarrow i+1$ and goes to repeat the above procedure from Step \ref{alg-steps: obtain-intermediate-estimate}.}
\end{algorithmic}
\end{algorithm}

In Algorithm \ref{alg: distributed_solution}, the maximum number of consultation hops $h$ controls the length of time consumed to diffuse information over the network in each iteration. By applying the procedure of Algorithm \ref{alg: distributed_solution}, each node is able to combine intermediate estimates from a subset of its $h$-hop neighborhood adaptively. This holds even if the parameter vector and the data and noise statistics change slowly in time. We evaluate the performance of this algorithm in Section \ref{subsec:simulation_arbitrary}.


\section{Simulation Results}\label{sec:Numerical-Example}

We illustrate the application of the mATC strategy to a tree and general graph network, and compare its performance with those of the non-cooperative, ATC diffusion and the centralized strategies. In the simulations, the step sizes and the discount factors are set as $\mu_{k}=0.08$ and $\nu_{k}=0.05$, respectively, for all $k\in\mathcal{N}$, and the numerical results are averaged over 1000 instances unless otherwise indicated.

\subsection{A simple tree network}\label{subsec: numerical-example-tree-network}
\begin{figure}
\begin{centering}
\includegraphics[scale=0.68]{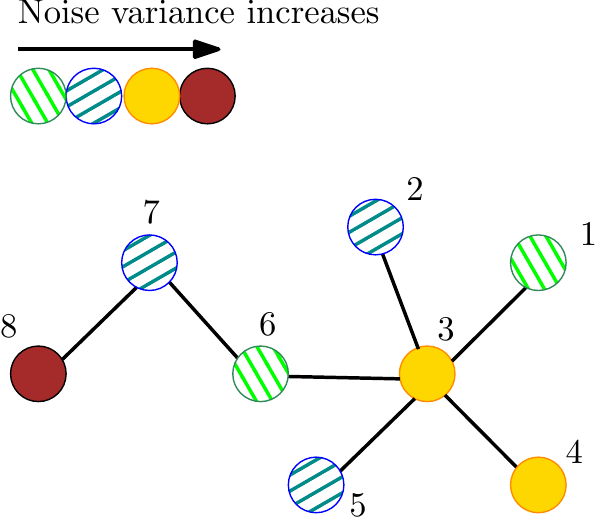}\,\,\,\,\,\,\,\,\includegraphics[scale=0.6]{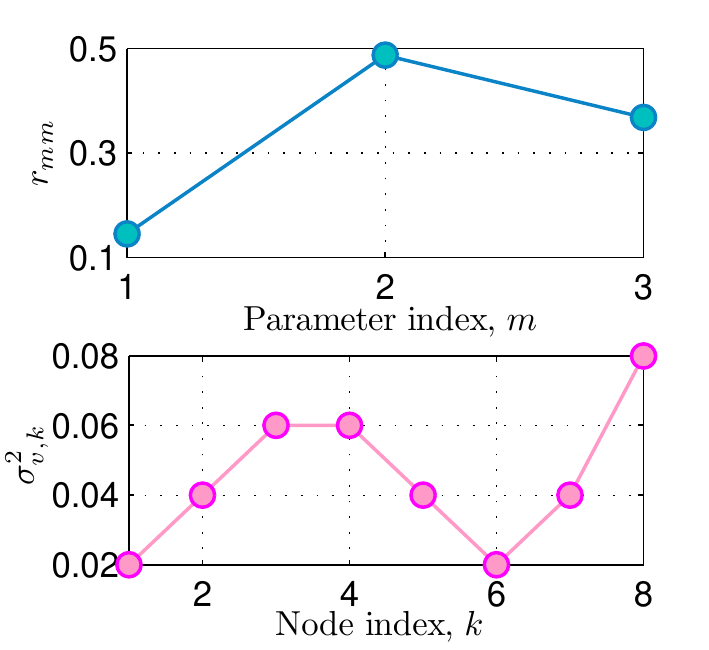}
\par\end{centering}
\caption{Network topology and data and noise profiles of every node. The number
next to a node denotes the node index.}
\label{fig:The-illustrative-example}
\end{figure}

An undirected tree network is shown in Fig.~\ref{fig:The-illustrative-example} together with the noise and data profiles of each node. The quantities $r_{mm}$ are
randomly generated such that $\sum_{m=1}^{M}r_{mm}=1$ and the regression
covariance matrix $R_{u,k}$ is diagonal with entries equal to $r_{mm}\cdot100\sigma_{v,k}^{2}$. Each node aims to estimate a $3\times1$ vector $\omega^{o}$ with every entry equal to $\nicefrac{1}{\sqrt{3}}$. We impose a maximum number $n_b$ of broadcasts allowed in the network per iteration as the network-wide energy constraint.

\begin{figure}
\begin{centering}
\includegraphics[scale=0.5]{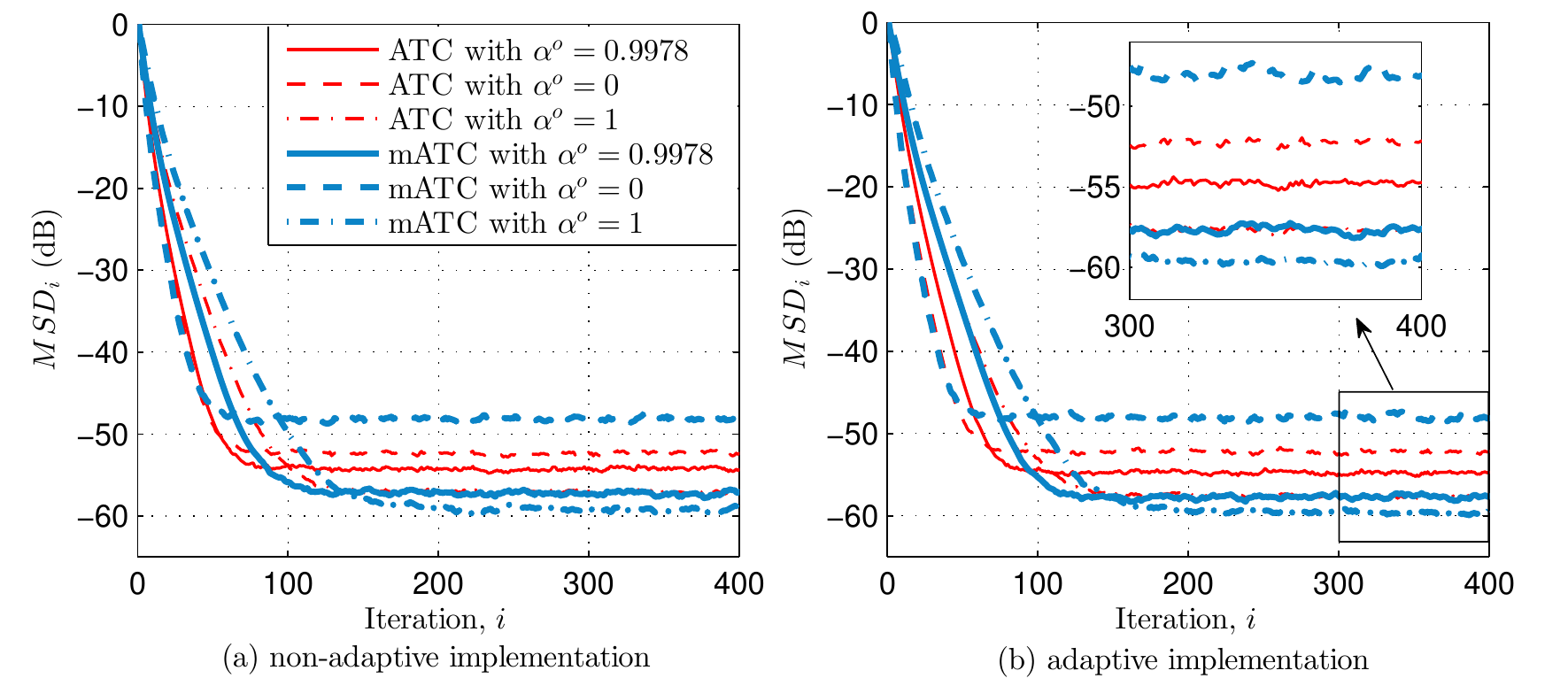}
\par\end{centering}
\caption{MSD performance of the ATC and the mATC diffusion strategies implementing the balancing rule of Theorem \ref{thm: Optimal-combination-rule}.}
\label{fig: The-illustrative-example: rule comparison}
\end{figure}

We illustrate the flexibility and usefulness of the proposed balancing rule in Theorem \ref{thm: Optimal-combination-rule} for assigning the combination weights. We set $n_{b}$ to 8, the same as that invoked by the ATC strategy. Both the ATC and the mATC strategies adopt the balancing rule
with the balancing coefficient $\alpha^o$ equal to the optimal value $0.9978$, or
the extreme value of 0 or 1. As shown in Fig.~\ref{fig: The-illustrative-example: rule comparison}, while the balancing rule with $\alpha^{o}=0$ gives the best transient network MSD, the rule with $\alpha^{o}=1$ gives the minimum steady-state network MSD. In comparison, the optimal balancing
rule using $\alpha^{o}=0.9978$ finds a good balance between
these two extremes. The results also show that the mATC strategy outperforms the ATC strategy in the steady state when $\alpha^{o}$ is optimal or equal
to 1.

\begin{figure}
\begin{centering}
\includegraphics[scale=0.55]{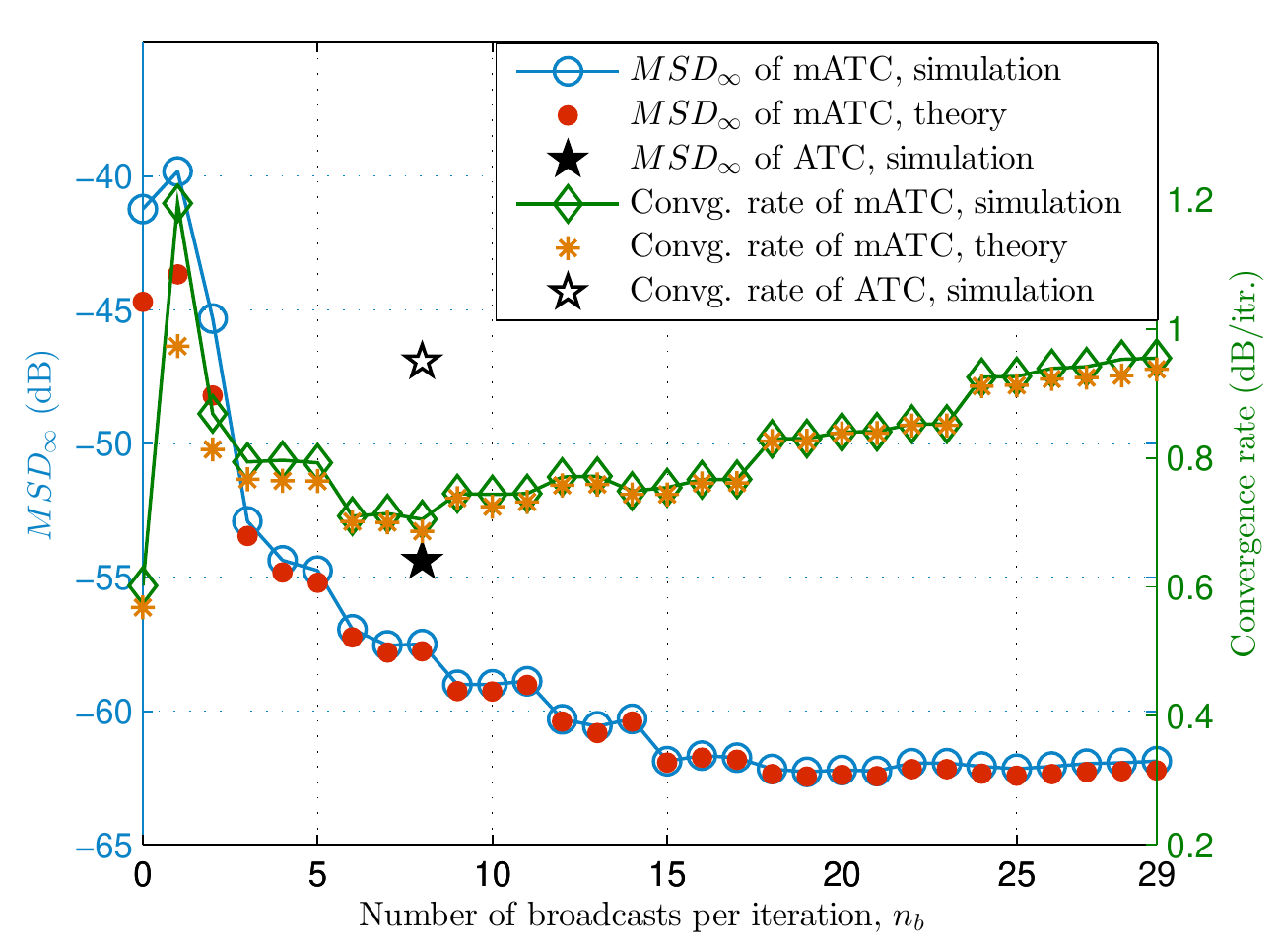}
\par\end{centering}
\caption{The steady-state value and the convergence rate of the network
MSD as functions of the number of broadcasts per iteration in the
network (where 0 and 29 broadcasts correspond to non-cooperative and centralized estimation, respectively).}
\label{fig:The-illustrative-example: trade-off curve}
\end{figure}

We next investigate the trade-off between the network performance
and the energy budget available per iteration. We define the \textit{convergence rate} as
the quotient of the decrease in the network MSD till 90\% of its
steady-state value divided by the number of iterations to achieve
that decrease. When the mATC adopts the optimal balancing rule, the theoretical (based on formula \eqref{eq: mean-square error evolution-2}) and numerical results are shown in Fig.~\ref{fig:The-illustrative-example: trade-off curve}. We observe a sharp increase in the convergence rate when the strategy transits
from a non-cooperative mode (with $n_{b}=0$) to the cooperative modes (with $n_{b}\ge1$). In the cooperative modes, the steady-state network MSD decreases almost monotonically as the number of broadcasts
increases, while the convergence rate first decreases and then increases. (The zig-zag variations in the curves are due to the limitations of minimizing the \textit{upper bound} of the steady-state network
MSD, and not a simulation artifact.) Once sufficient number of broadcasts are invoked, the marginal benefit brought by more consultations is small to the steady-state network
MSD, but still notable to the convergence rate.

\begin{figure}
\begin{centering}
\includegraphics[scale=0.58]{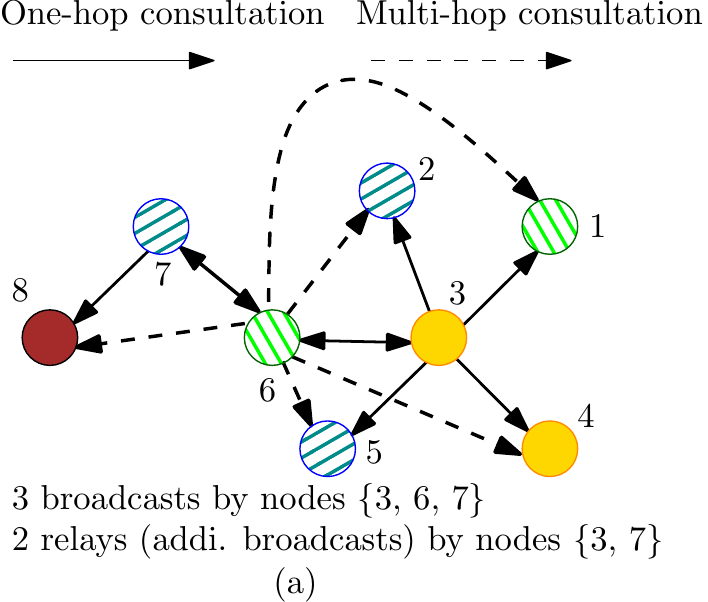}\includegraphics[scale=0.5]{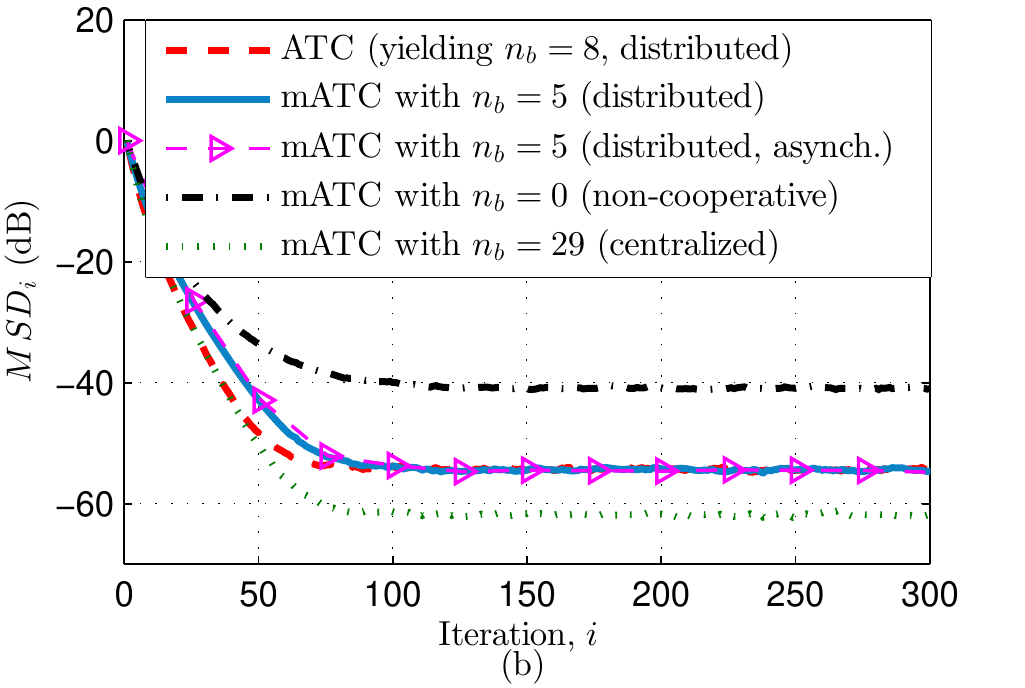}
\par\end{centering}
\noindent \centering{}\caption{Comparison of the ATC and mATC diffusion strategies adopting the optimal balancing
rule: (a) consultations in each iteration with the mATC strategy limited
to 5 broadcasts per iteration; (b) network MSD curves obtained
from simulations. \label{fig: The-illustrative-example: n5}}
\end{figure}

A specific solution of the information neighbor configuration is shown in Fig.~\ref{fig: The-illustrative-example: n5}, which invokes 3 or 37.5\% less broadcasts in each iteration compared to the ATC strategy. The figure also shows results obtained with an \emph{asynchronous} mATC strategy, in which intermediate estimates from neighbors more than one hop away are not relayed within the same diffusion iteration to a node. Instead, such estimates may arrive at the node in a later iteration, and combined only at that iteration. The asynchronous strategy avoids waiting times for receiving the two-hop information, and allows the mATC strategy to perform each combination step in the same time scale as the ATC strategy. The estimation performance turns out to be close (i.e., the difference is almost negligible) to the synchronous mATC strategy.

\subsection{An arbitrary network}\label{subsec:simulation_arbitrary}

\begin{figure}
\begin{centering}
\includegraphics[scale=0.51]{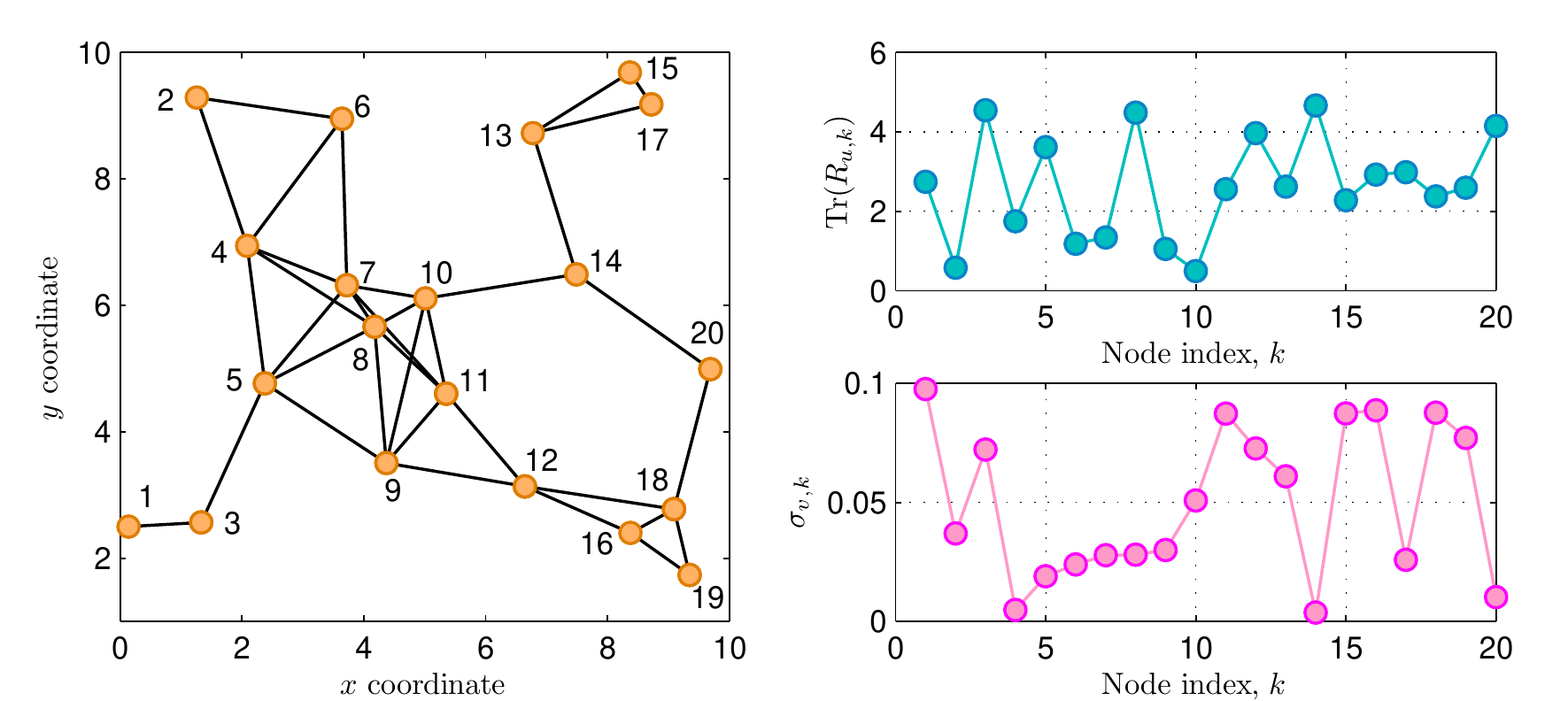}
\par\end{centering}
\caption{The topology of a random network, and the data and noise profiles
at each node. The number next to a node is the node index.\label{fig:The-random-example}}
\end{figure}

We randomly generated an undirected network with 20
nodes within a $10\times10$ square area, as shown in Fig.\
\ref{fig:The-random-example} together with the noise and data power profiles, where the data power $\mathrm{Tr}(R_{u,k})$
is equally distributed over the parameter components of each node $k$. Each node aims to estimate a $2\times1$ vector $\omega^{o}$ with every entry equal to $\nicefrac{1}{\sqrt{2}}$. To mimic a heterogenous environment, we assume that the communication cost incurred is proportional
to the square distance of a node to its farthest directly reachable neighbor, i.e., the communication cost coefficient of node $k$ in \eqref{eqn:commcost} is given by
\[
c_{k}^{cm,0}=c^{0}\times\left(\max\{\text{dist}(k,l) : l\in\dreachN{k}\}\right)^{2},
\]
where $c^{0}$ is a given scalar and $\text{dist}(k,l)$ is the distance between node $k$ and node $l$. Without loss of generality, we set $c^{0}=1$, and we study the energy-performance trade-off when the network is subject to a network-wide energy budget.

\begin{figure}
\begin{centering}
\includegraphics[scale=0.55]{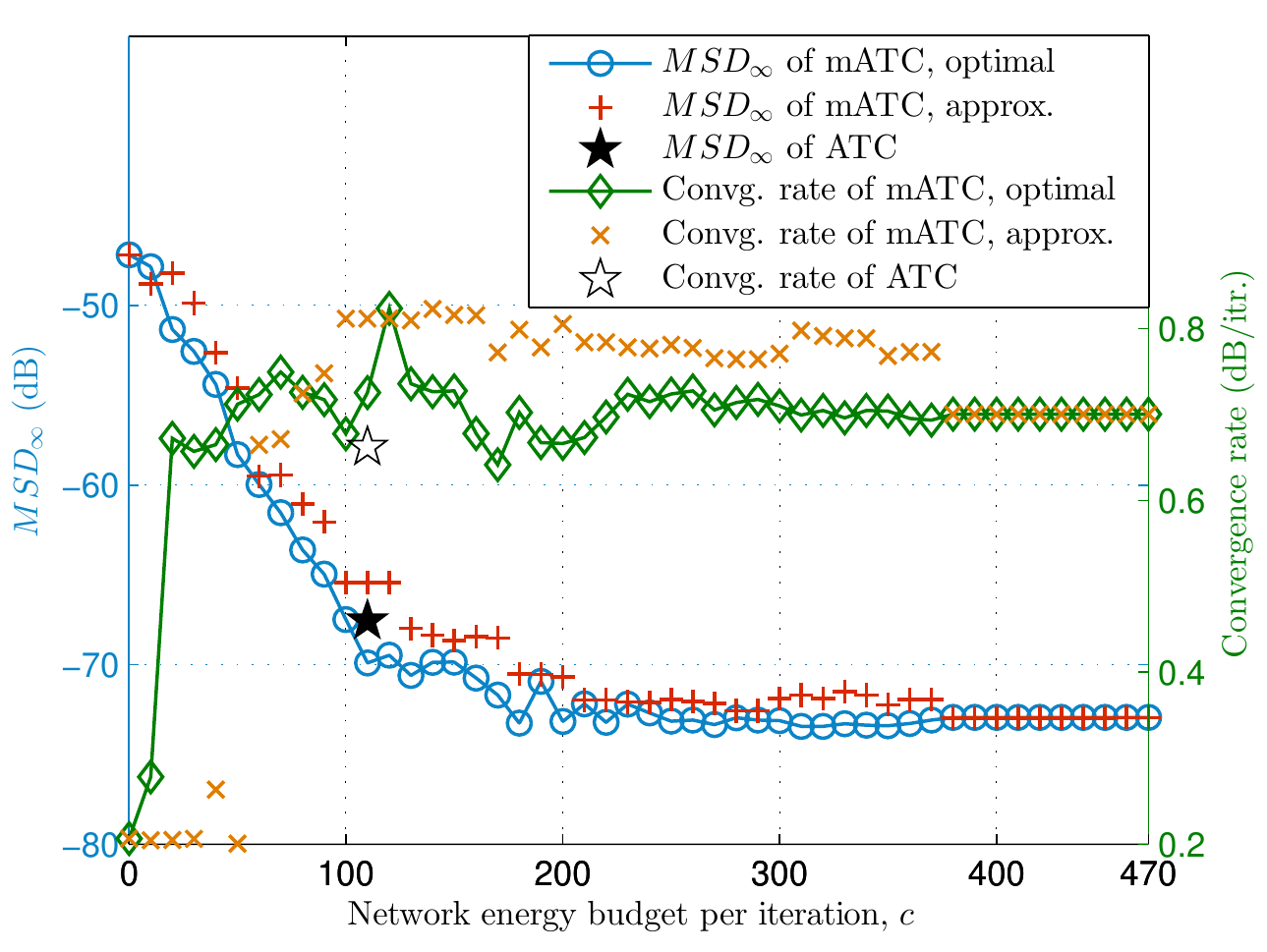}
\par\end{centering}
\caption{The energy-performance trade-offs and convergence rates of mATC and ATC strategies.}
\label{fig:The-random-example: trade-off curve}
\end{figure}

\begin{figure}
\centering
        \begin{subfigure}[tbp]{0.24\textwidth}
                \includegraphics[width=\textwidth]{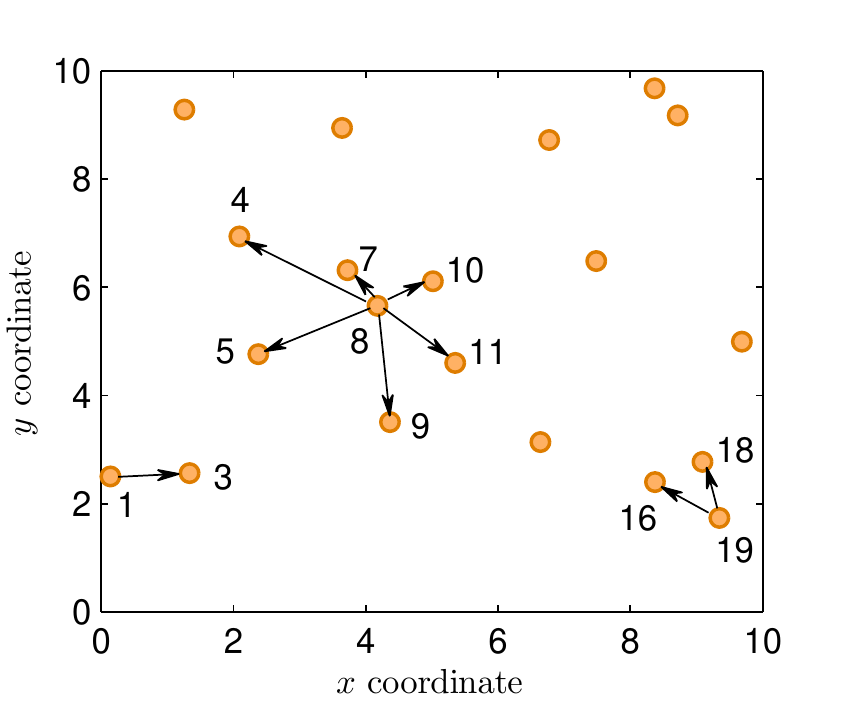}
                \caption{Optimal info neighbor configuration with $c=10$.}
                \label{fig:c10}
        \end{subfigure}
        \begin{subfigure}[tbp]{0.24\textwidth}
                \includegraphics[width=\textwidth]{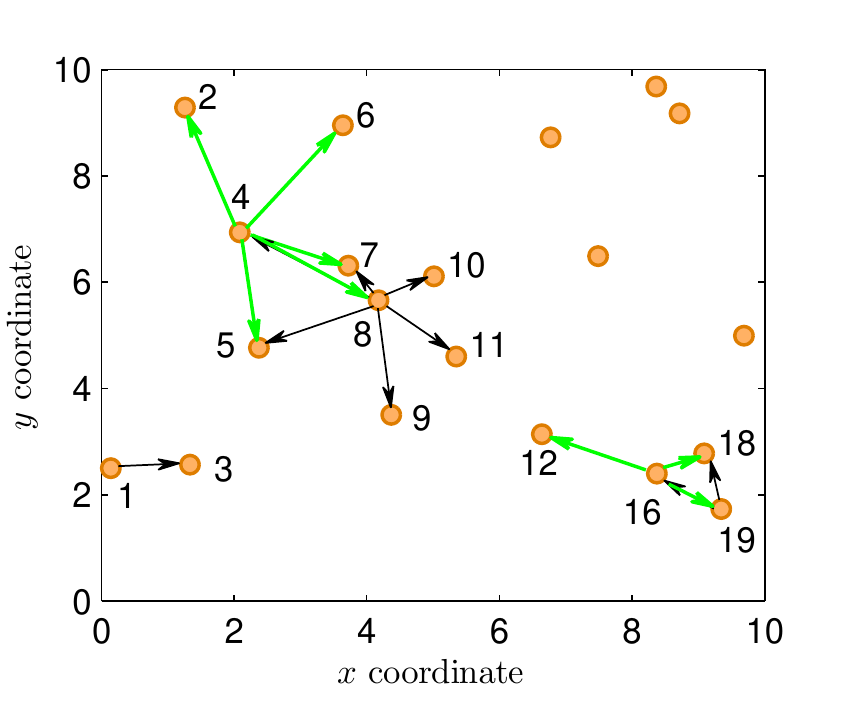}
                \caption{Optimal info neighbor configuration with $c=20$.}
                \label{fig:c20}
        \end{subfigure}

        \begin{subfigure}[tbp]{0.24\textwidth}
                \includegraphics[width=\textwidth]{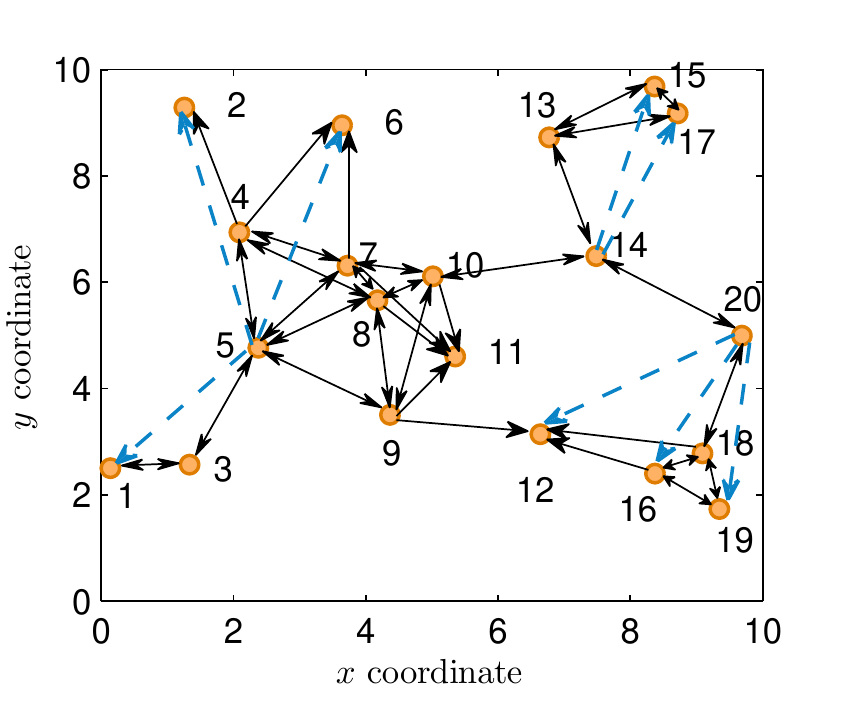}
                \caption{Optimal info neighbor configuration with $c=110$.}
                \label{fig:c110}
        \end{subfigure}
        \begin{subfigure}[tbp]{0.24\textwidth}
                \includegraphics[width=\textwidth]{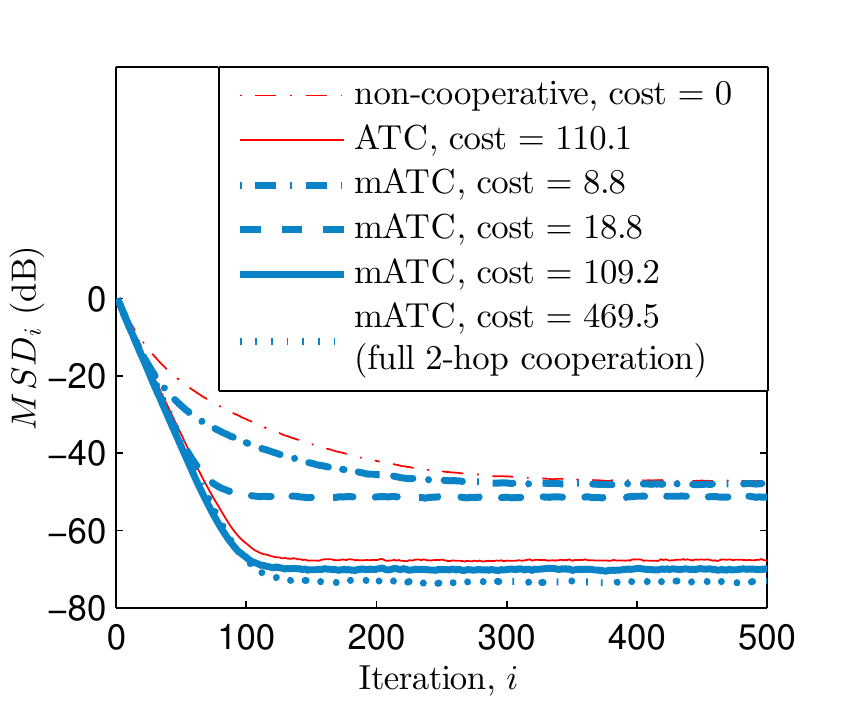}
                \caption{Network MSD curves for different energy costs.}
                \label{fig:MSD for c10-20-110}
        \end{subfigure}
\caption{Three optimal configurations of information neighbors and related network MSD curves obtained from simulations. The arrows in (a)--(c) indicate the diffusion directions. The thicker green arrows in (b) indicate new diffusions relative to those in (a), and the light blue dashed arrows in (c) represent two-hop diffusions that require neighbors' relay.}
\label{fig: The-random-example: three instances}
\end{figure}

Fig.~\ref{fig:The-random-example: trade-off curve} shows the energy-performance trade-off when the mATC diffusion strategy adopts the balancing rule of Theorem \ref{thm: Optimal-combination-rule}. We show the performances of both the exact solution found using the MILP (P3) and the approximate solution found through relaxations using Algorithm \ref{alg: appro_solution}. We observe that the steady-sate network MSD decreases almost monotonically with the energy budget, while the convergence rate of the network MSD first increases as cooperation is enabled, and then fluctuates and becomes steady as more energy budget is available. The corresponding changes in the information neighbor configuration are illustrated in Fig.~\ref{fig:c10}--\ref{fig:c110}. The performance gain turns out to be very small if the energy budget is large enough (larger than 250 in this case). The approximate solutions give similar energy-performance trade-offs, but have uniformly worse steady-state network MSDs.

In addition to ATC, we also compare the mATC strategy with the game theoretic combine-then-adapt diffusion strategy proposed in \cite{namvar2013distributed}. We call this the CTA-game strategy for short. This strategy requires the setting of a couple of parameters, and we refer the reader to \cite{namvar2013distributed} for the full details. Specifically, the local and global utility gains $K_{l,1}$ and $K_{g}$ are chosen to be 1, the local energy price $K_{l,2}$ is given a value in the range of [0, 0.5] for every sensor (the network steady-state MSD and the convergence rate are found to remain the same once $K_{l,2}\ge 0.5$ in this case), the ``exploration'' factor $\delta$ is set to 0.1, and the other parameters $\gamma_{l}$ and $\gamma_{g}$ are fixed as 1. For the CTA-game strategy, we average the simulations results over 300 random instances, each using 2000 iterations.

\begin{figure}
\begin{centering}
\includegraphics[scale=0.6]{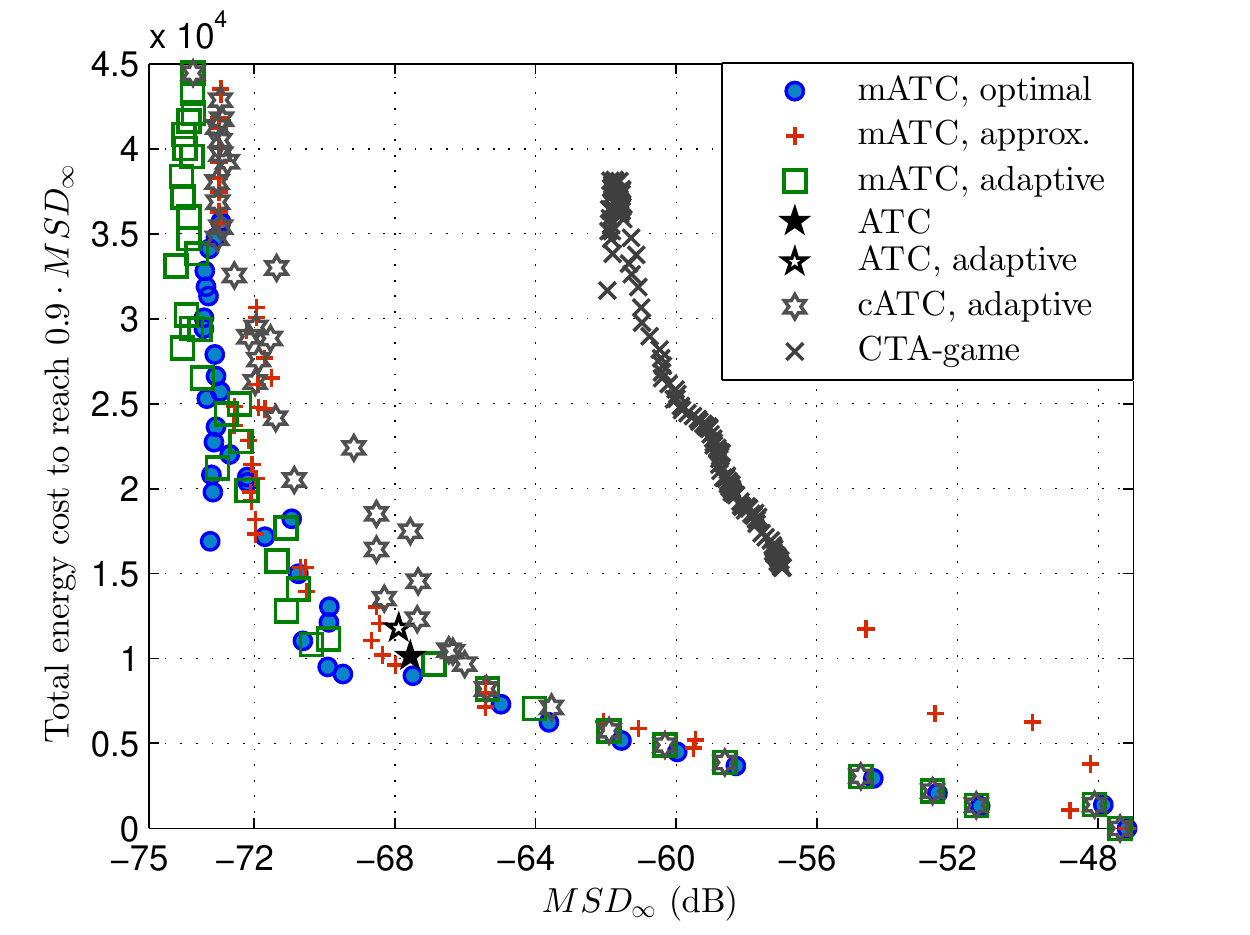}
\par\end{centering}
\caption{Average total energy cost required for the network to converge to $90\%$ of its steady-state MSD value.}
\label{fig:The-random-example: total_energy_comp}
\end{figure}

In Fig. \ref{fig:The-random-example: total_energy_comp}, we compare the CTA-game strategy with two versions of the mATC strategy: one implemented using full knowledge of the noise and input data statistics, and another in an adaptive manner per Algorithm \ref{alg: distributed_solution} (where the local energy budgets correspond to the configurations that yield the optimal performance-energy trade-off curve shown in Fig. \ref{fig:The-random-example: trade-off curve}). We also compare against an adaptive implementation of an energy-constrained ATC strategy, which we call cATC, in which each node diffuses only its own intermediate estimate if its energy budget permits. We see that both versions of mATC require much less total energy to converge to 90\% of the same steady-state network MSD as the CTA-game strategy. For example, for a steady-state MSD of -62 dB, the adaptive mATC strategy consumes about $0.6\times 10^4$ units of energy. In contrast, the CTA-game strategy consumes about $3.6\times 10^4$ units of energy, which is six times of that used by the adaptive mATC. The higher total energy consumption of the CTA-game strategy is found to be caused by its slow convergence, which does not change much even if the parameters are tuned in the given ranges. The simulation results also indicate that the adaptive mATC diffusion strategy requires less energy for the network to converge to the same MSD, as compared to the adaptive cATC strategy.

\begin{figure}
\begin{centering}
\includegraphics[scale=0.6]{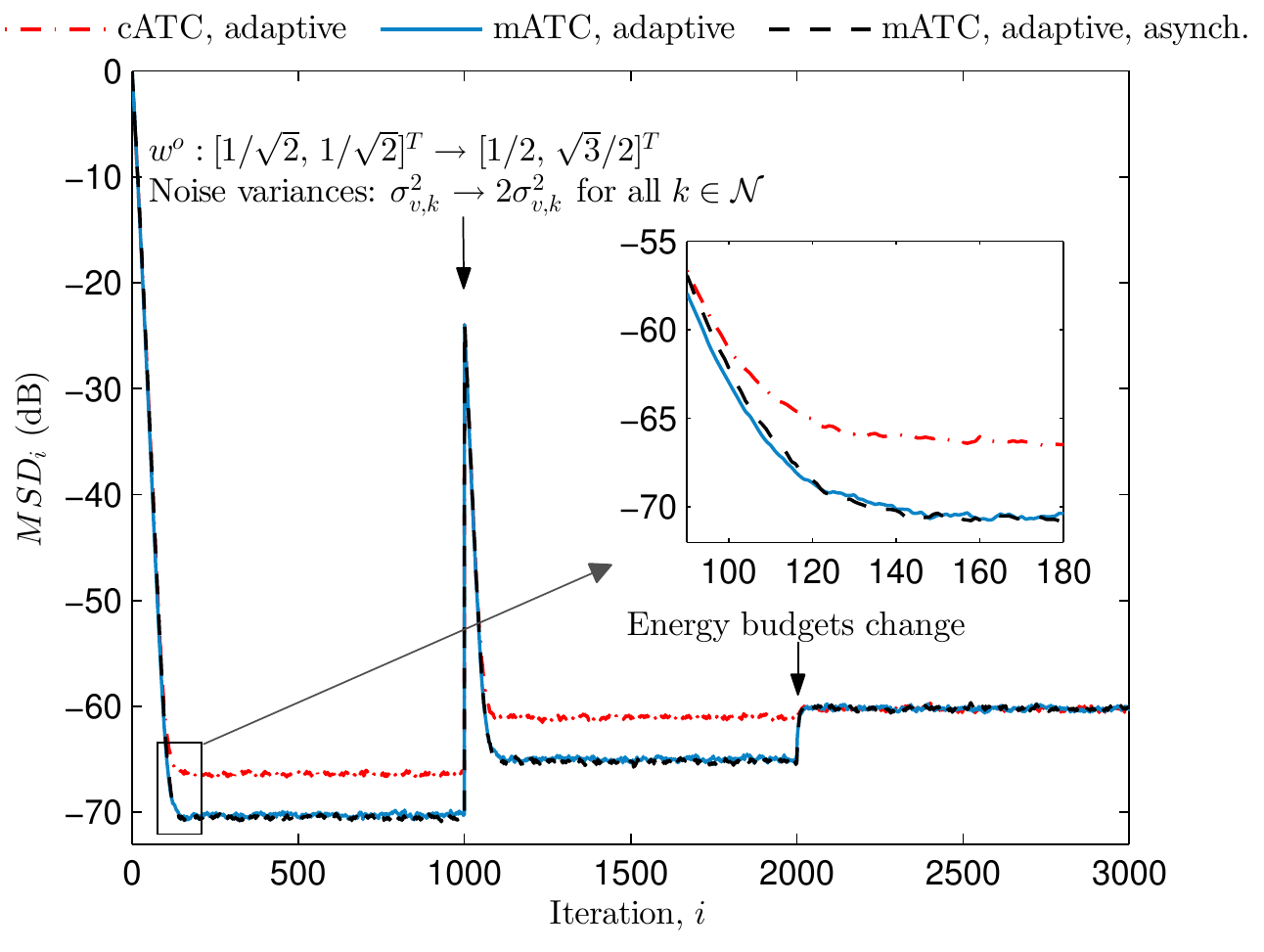}
\par\end{centering}
\caption{Performances of the cATC and mATC diffusion strategies in response to changes in network conditions.}
\label{fig:The-random-example: adaptability of mATC}
\end{figure}

We also examine the adaptability of the mATC strategy in response to changes in the parameter $w^o$, the measurement noise variances and the locally available energy budgets. We assume that the network starts with the configuration shown in Fig.~\ref{fig:c110}, where each node has the exact local energy budget to support the transmissions shown. We implement an asynchronous version of the adaptive mATC strategy in Algorithm~\ref{alg: distributed_solution} (i.e., intermediate estimates from multi-hop neighbors are not relayed within one diffusion iteration), and the adaptive cATC strategy. After 1000 iterations, the parameter vector $w^o$ changes from the original $[1/\sqrt{2}, 1/\sqrt{2}]^T$ to the new $[1/2, \sqrt{3}/2]^T$ while the measurement noise variance $\sigma_{v,k}^2$ is doubled for all $k\in\mathcal{N}$. After 2000 iterations, the local energy budgets available to the nodes 3, 4, 13 and 18 all decrease to support only one broadcast instead of two broadcasts per iteration, and the local energy budgets available to nodes 7 and 8 are set to zero. We observe from Fig.~\ref{fig:The-random-example: adaptability of mATC} that ATC, mATC, and the asynchronous mATC are all able to adapt to the changes in network conditions, and that both mATC and asynchronous mATC achieves steady-state MSDs better than cATC (except after the last change when mATC becomes equivalent to ATC), thanks to their selective multi-hop consultations that fully exploit the local energy budgets.

\section{Conclusion\label{sec: Conclusion}}

In this paper, we have considered the use of multi-hop diffusion that allows nodes to exchange intermediate parameter estimates with their selected information neighbors instead of just the physical neighbors. For two classes of networks, we propose an MILP to select the information neighbors together with the relay nodes for each node, which approximately optimizes a trade-off between the available energy budgets for each iteration, and the steady-state network MSD performance. For arbitrary networks in which there are only local energy budget constraints, and consultations constrained to within a fixed number of hops, we propose a distributed and adaptive algorithm to select the information neighbors. Simulation results suggest that our proposed methods achieve better MSD performance than the traditional diffusion strategy, while having a same or lower communication cost.

Our current optimization procedure for networks with a network-wide energy budget requires knowledge of the network topology as well as the data and noise variance profiles of every node. This implies that the optimization can only be performed at a centralized processor, and only infrequently. It would be of future research interest to develop distributed optimization techniques like those in \cite{towfic2014adaptive} to perform online adaptive optimization as the network conditions vary over time, and to study the frequency at which such optimization needs to be run, in order to maintain a reasonable level of optimality.

\section*{Appendix A}

\section*{Derivation of the SDP in \eqref{eq: SDP for beta} for solving $\beta$}

Some well-known matrix results used by the derivation are first given in the
following lemma, in which (a) is proved by using the singular value decomposition technique and the proofs of (c)-(d) can be found in pages 399, 19 and 473 of \cite{Horn1990}, respectively.

\begin{lem}\label{lemma:elem_matrix}\
\begin{enumerate}[(a)]
	\item\label{it:XX} For any matrix $X\in\mathbb{R}^{n\times n}$, $X^{T}X\preceq I_{n}$
if and only if $XX^{T}\preceq I_{n}$.
	\item\label{it:CAC} If $A\in\mathbb{R}^{n\times n}$ is positive semi-definite, then for any $C\in\mathbb{R}^{n\times n}$, $C^{*}AC$ is positive semi-definite.
	\item\label{it:XY} If matrices $X$, $Y$ and $X+Y$ are invertible, then $(X+Y)^{-1}=X^{-1}-X^{-1}(Y^{-1}+X^{-1})^{-1}X^{-1}$.
	\item\label{it:CBAB} Suppose that a Hermitian matrix is partitioned as
$\begin{smallmatrix}\left[\begin{array}{cc}
A & B\\
B^{*} & C
\end{array}\right]\end{smallmatrix}$, where $A$ and $C$ are positive definite. This matrix is positive semi-definite if and only if $C-B^{*}A^{-1}B$ is positive semi-definite.
\end{enumerate}
\end{lem}

{
Let $P_\beta\triangleq\mathcal{MSM}/\beta$ and $Q\triangleq(I_{NM}-\mathcal{MR})^{2}$,
which are positive definite matrices; and let $\phi_\beta\triangleq\sqrt{P_\beta+Q}$, which is a positive
definite matrix, and so $\phi_\beta^{2}= P_\beta+Q$.
}
The derivation proceeds
as follows (other constraints on the combination weight matrix $A$
and the constraint $\beta>0$ are not shown):
\begin{align*}
 & \min_{A,\,\beta}\beta,\,\,
\text{s.t. }\lambda_{\max}\left(\mathcal{Y}/\beta+\mathcal{B}\mathcal{B}^{*}\right)\le1,\,\, A^{T}\ones{N}=\ones{N}\\
 & \Leftrightarrow\min_{A,\,\beta}\beta,\,\,\text{s.t. }\frac{\mathcal{Y}}{\beta}+\mathcal{B}\mathcal{B}^{*}\preceq I_{NM},\,\, A^{T}\ones{N}=\ones{N}\\
 & \Leftrightarrow\min_{A,\,\beta}\beta,\,\,\text{s.t. }\mathcal{A}^{T}\phi_\beta^{2}\mathcal{A}\preceq I_{NM},\,\,,\,\, A^{T}\ones{N}=\ones{N}\\
 & \overset{\text{Lemma \ref{lemma:elem_matrix}\eqref{it:XX}}}{\Longleftrightarrow}\min_{A,\,\beta}\beta,\,\,\text{s.t. }\phi_\beta\mathcal{A}\mathcal{A}^{T}\phi_\beta\preceq I_{NM},\,\, A^{T}\ones{N}=\ones{N}\\
 & \overset{\text{Lemma \ref{lemma:elem_matrix}\eqref{it:CAC}}}{\Longleftrightarrow}\min_{A,\,\beta}\beta,\,\,\text{s.t. }\mathcal{A}\mathcal{A}^{T}\preceq\phi_\beta^{-2},\,\, A^{T}\ones{N}=\ones{N}
\end{align*}
By left multiplying both sides of the SDP constraint with $\ones{N}^{T}\otimes I_{M}$
and right multiplying with $\ones{N}\otimes I_{M}$,
the last optimization implies that an approximate solution of $\beta$ can
be solved from:{\small
\begin{align*}
& \min_{\beta}\beta,\,\,\text{s.t. }(\ones{N}^{T}\otimes I_{M})(\ones{N}\otimes I_{M})\\
& \quad\quad\quad \preceq(\ones{N}^{T}\otimes I_{M})(Q+P_\beta)^{-1}(\ones{N}\otimes I_{M})\\
 & \overset{\text{Lemma \ref{lemma:elem_matrix}\eqref{it:XY}}}{\Longleftrightarrow}\min_{\beta}\beta,\,\,\text{s.t. }(\ones{N}^{T}\otimes I_{M})(\ones{N}\otimes I_{M})\\
 & \,\,\,\,\,\,\preceq(\ones{N}^{T}\otimes I_{M})[Q^{-1}-Q^{-1}(P_\beta^{-1}+Q^{-1})^{-1}Q^{-1}](\ones{N}\otimes I_{M})\\
 & \Leftrightarrow\min_{\beta}\beta,\,\,\text{s.t. }(\ones{N}^{T}\otimes I_{M})(Q^{-1}-I_{NM})(\ones{N}\otimes I_{M})\\
 & \,\,\,\,\,\,-(\ones{N}^{T}\otimes I_{M})Q^{-1}(P_\beta^{-1}+Q^{-1})^{-1}Q^{-1}(\ones{N}\otimes I_{M})\succeq0\\
 & \overset{\text{Lemma \ref{lemma:elem_matrix}\eqref{it:CBAB}}}{\Longleftrightarrow}\min_{\beta}\beta,\,\,\text{s.t. }\\
 & \left[\begin{array}{cc}
P_\beta^{-1}+Q^{-1} & Q^{-1}(\ones{N}\otimes I_{M})\\
(\ones{N}^{T}\otimes I_{M})Q^{-1} & (\ones{N}^{T}\otimes I_{M})(Q^{-1}-I_{NM})(\ones{N}\otimes I_{M})
\end{array}\right]\succeq0.
\end{align*}
}The last optimization problem is an SDP, and can be solved by standard solvers to get an approximate solution of $\beta$ for the original problem.

\section*{Proof of Theorem \ref{thm: validity of Algorithm 1}}

We first prove the theorem for (P2). The proof relies on the following lemma, which can be shown by using the valid inequalities \eqref{eq: vc-pi-pi-P2}, and is omitted for brevity.

\begin{lem}\label{lm: Order-relationship-relay-variables}
Let the optimal relay and selection variables of the relaxed LP corresponding to (P2) (including the valid inequalities in (\ref{eq: vc-pi-pi-P2})) be $\{(\tilde{\pi}_{lk},\tilde{\delta}_{lk}): l\in\multiN{k}, k\in\cN \}$. For any relay path from a node $l$ to another node $j$ such that $\tilde{\delta}_{lj}=1$, and any pair of nodes $(k_{1},\, k_{2})$, where $k_{2}$ is a successor of $k_{1}$ on the relay path, we have $\tilde{\pi}_{lk_{1}}\geq\tilde{\pi}_{lk_{2}}$.
\end{lem}

To prove that the solution of the relay and selection variables returned by Algorithm \ref{alg: appro_solution} is feasible for (P2), we need to show that: (i) the relay variables, obtained after thresholding those from the relaxed LP, result in valid transmission paths (i.e., every path is able to deliver the information as desired); and (ii) the local and network-wide energy budgets are satisfied
by each node and the whole network, respectively. These are proved as follows.

Observe that line \ref{alg-steps: initial-regularized-solution} of Algorithm \ref{alg: appro_solution} returns relay variables feasible for (P2) without energy budget constraints. This is because the solution is able to indicate every information relay path without ambiguity as shown by Lemma \ref{lm: Order-relationship-relay-variables}.

Lines \ref{alg-steps: global-budget-satisfication-start}--\ref{alg-steps: global-budget-satisfication-end} choose relays in non-decreasing order of their corresponding relay variable values, and removes them so as to satisfy the network-wide energy budget constraint. By Lemma \ref{lm: Order-relationship-relay-variables}, the operation only removes tails of certain relay paths sequentially until the network-wide energy budget constraint is satisfied. Therefore, the remaining part of the relay path is still valid.

In lines \ref{alg-steps: local-budget-satisfication-start}--\ref{alg-steps: local-budget-satisfication-end}, selected relays for each node $k\in\cN$ are removed in non-decreasing order of their $\tilde{\pi}_{lk}$ values in order to reduce node $k$'s energy consumption to not more than $c_k$ so that \eqref{eq: cons-node-wide-energy} is satisfied. If node $k$ stops relaying information from node $l^*$, then we also set $\pi_{l^*j}=0$ for all nodes $j$ that used to obtain node $l^*$'s information through node $k$. This ensures that the relay path remains valid.

Since a feasible solution of the relay variables uniquely determines a feasible solution of the selection variables, the solution of the selection variables obtained from line \ref{lm: Order-relationship-relay-variables} is feasible for (P2).

The theorem for (P3) can be proved in a similar way, because the result in Lemma \ref{lm: Order-relationship-relay-variables} remains true for the relaxed (P3) with the valid inequalities \eqref{eq: vc-delta-implied-by-pi-P3}. The proof is now complete.

\bibliographystyle{IEEEtran}
\bibliography{HuTay_gATC_2014}

\end{document}